\documentclass[12pt,twoside]{article}

\usepackage{amsmath,amsfonts,amssymb,theorem}

\newtheorem{ques}{Question}
\newtheorem{theo}{Theorem}
\newtheorem{prop}{Proposition}
\newtheorem{coro}{Corollary}
\newtheorem{lem}{Lemma}

\theorembodyfont{\rmfamily}

\newenvironment{prv}{\paragraph{Proof}}{\par\medskip}

\begin{document}

\author{Philippe Eyssidieux \\{\small CNRS - Laboratoire Emile Picard, Universit\'e 
Paul Sabatier, Toulouse, France} \\ {\footnotesize eyssi@picard.ups-tlse.fr}
\and  Ngaiming Mok \footnote{Research 
supported by a CERG grant of the Research Grants Council of Hong-Kong} \\
{\small Dept. 
of Math., Hong Kong University, Hong-Kong, China}\\ {\footnotesize 
nmok@hkucc.hku.hk}}

\date{\today}

\title{On the Validity or Failure of Gap Rigidity \\
  for Certain pairs of \\
  Bounded Symmetric Domains}

\maketitle

\begin{flushright}
{\it in memory of Professor Armand Borel}
\end{flushright}

Let $\Omega$ be a bounded symmetric domain equipped with a canonical 
K\"ahler metric.  We are interested to characterize 
holomorphic geodesic cycles (i.e., compact complex geodesic submanifolds)
$S \subset X$ on Hermitian locally
symmetric manifolds $X$ uniformized by $\Omega$, i.e., $X = \Omega/\Gamma$, 
where
$\Gamma \subset \text{Aut}(\Omega)$ is any  torsion-free discrete group of 
automorphisms,
in terms of differential-geometric or algebro-geometric conditions.  In 
Eyssidieux-Mok \cite{EM} we studied almost geodesic complex submanifolds and 
formulated the gap phenomenon.  
Up to equivalence under $\text{Aut}(\Omega)$ there 
are only a finite number of totally geodesic complex submanifolds $D \subset 
\Omega$ which
are themselves biholomorpic to bounded symmetric domains.  Let $\epsilon > 0$.
We say that $S \subset X$ is $\epsilon$-geodesic if and only if the norm of 
the 
second fundamental form of $S \subset X$ is uniformly bounded by $\epsilon$.  
Fixing
$\Omega$ but letting $\Gamma \subset \Omega$ be arbitrary we showed that when 
$\epsilon$ is sufficiently small, an $\epsilon$-geodesic compact complex
submanifold $S \subset X$ must locally 
resemble 
one 
and 
only 
one embedding
of 
the 
models $(\Omega,D;i)$, $i: D \hookrightarrow \Omega$.  We say that gap rigidity
holds for the pair $(\Omega,D;i)$ if for $\epsilon$ sufficiently small,
any $\epsilon$-geodesic compact complex submanifold $S \subset X$ modeled
on $(\Omega,D;\imath)$ is necessarily totally geodesic for any $X$ uniformized
by $\Omega$.  When $i: D \hookrightarrow \Omega$ is understood, we will 
sometimes just write $(\Omega,D)$.

\noindent There is a strong form of the gap phenomenon, where one can 
characterize
certain
holomorphic geodesic cycles in algebro-geometric terms, in terms of the 
genericity of their tangent spaces.  Write $\Omega = G/K$, where $G = 
\text{Aut}(\Omega)$
and $K \subset G$ is the isotropy subgroup at a reference point $o \in \Omega$.
We will say that gap rigidity holds for $(\Omega,D)$ in the Zariski topology
if there 
exists a $K$-invariant Zariski-open subset $\mathcal 
O_o$ 
in 
the 
Grassmannian 
$\text{{\rm Gr}}(\text{dim}(D),T_o(\Omega))$ of $\dim(D)$-planes in 
$T_o(\Omega)$ 
such that 
$[T_o(D)] \in {\mathcal O}_o$, and such that, for any 
complex 
manifold 
$X$ 
uniformized
by $\Omega$, any compact complex submanifold $S \subset X$ of dimension equal 
to 
dim$(D)$ must be totally geodesic, whenever every tangent plane $[T_x(S)]$ of 
$S$
lifts to an element of $\mathcal O_o$.  
Obviously 
if 
gap 
rigidity 
holds 
for 
$(\Omega,D)$
in the Zariski topology, it holds for $(\Omega,D)$ in the 
differential-geometric
sense as explained.  To make a distinction we will sometimes refer to the 
latter as
gap rigidity in the complex topology.

\noindent The simplest example of gap rigidity in the Zariski topology is the 
case
of a product domain $(D^k,D;\delta)$, where $\delta$ always refers to the 
diagonal 
embedding, given by $\delta(z) = (z,\cdots,z)$, as given in \cite{EM}.  The 
observation there resulted from 
the uniqueness of K\"ahler-Einstein metrics on a compact complex manifold
with ample canonical line bundle.  Eyssidieux 
\cite{Eys1}\cite{Eys2} 
considered, 
in the context
of variations of Hodge strutures, the question of characterizing certain 
holomorphic
geodesic cycles by Chern number inequalities, establishing as a by-product gap 
rigidity in the Zariski topology for many pairs of $(\Omega,D)$ by means of 
Gauss-Manin complexes. His methods apply to the case of period domains arising 
from Hodge theory.  
In Mok \cite{M2} using intersection theory on the projectivized tangent bundle 
we 
determined the set of all pairs 
$(\Omega,D)$ with $\Omega$ irreducible and $\dim(D)=1$  for which gap rigidity 
holds in the 
Zariski topology.  They are precisely those $\Omega$ of rank $r > 1$ such that 
the highest characteristic subvariety (cf. \cite{M2},\S1) 
$\mathcal S_o \subset \Bbb P(T_o(\Omega))$ is a hypersurface in $\Bbb 
P(T_o(S))$.
At
the same time, exploiting existence and uniqueness of the K\"ahler-Einstein 
manifolds with ample canonical line 
bundle, 
G-structures,
and 
Hermitian
metric rigidity, we proved that for $\Omega = D^{IV}_n$ the $n$-dimensional 
bounded symmetric domain of type IV (dual to the hyperquadric) and of 
dimension $\geq 3$, gap rigidity holds in the Zariski topology for 
$(D^{IV}_n,D^{IV}_k,i)$, where $k \geq 2$ and $i:D^{IV}_k\hookrightarrow
D^{IV}_n$ is the standard embedding.

\noindent The purpose of the article is two-fold.  First of all, we will show 
that in
general 
gap rigidity already fails in the complex topology.  More precisely, we show
that gap rigidity fails for $(\Delta^2,\Delta\times \{0\})$ by constructing 
a sequence of ramified coverings $f_i: S_i \to T_i$ between hyperbolic compact
Riemann surfaces such that, with respect to norms defined by the Poincar\'e
metrics, sup$||df_i|| \to 0$ as $i \to \infty$.  Since any bounded symmetric
domain of rank $\geq 2$ contains a totally geodesic bidisk, this implies
that gap rigidity fails in general on any bounded symmetric domain of 
rank $\geq 2$.  Our counterexamples make it all the more interesting to find
sufficient conditions for pairs $(\Omega,D)$ for which gap rigidity holds. This
will be addressed in the second part of the article, where for $\Omega$
irreducible, we generalize the results for holomorphic curves in \cite{M2} to 
give
a sufficient condition for gap rigidity to hold for $(\Omega,D)$ in the 
Zariski 
topology. In notations of the preceding paragraphs, we show that it is 
sufficient    
to take $\mathcal O_o \subset T_o(\Omega)$ to be such that the 
excluded
subvariety
${\mathcal Z}_o = \text{{\rm Gr}}({\dim(D)} ,T_o(\Omega)) -
{\mathcal
O}_o$ 
is a
hypersurface.  In the 
terminology of geometric invariant theory we say  
that $[T_o(D)] \in \text{{\rm Gr}}({\text{dim}}(D),T_o(\Omega))$ is a
semistable
point.
For the case of bounded symmetric domains, this criterion applies to the case 
of
1-hyperrigid domains $D \subset \Omega$ of \cite{Eys1}\cite{Eys2} and to the
examples of
\cite{M2} for holomorphic curves and from quadric structures, giving a 
unifying 
conceptual explanation and a more general framework for the validity of the 
phenomenon. We will show that the class of GIT-semistable embeddings of 
irreducible bounded symmetric domains
coincides with the class of $(H_2)$-embeddings, see \cite{Sat}. 
In \cite{Sat1}, I. Satake classified $(H_2)$-embeddings into classical domains 
and
the full classification has been obtained 
 by S. Ihara \cite{Iha} as a step towards his
classification of all embeddings of bounded symmetric domains.
Since the abstract classification of GIT-semistable embeddings does not give an
explicit invariant hypersurface in each case, we will also give a detailled 
exposition
of some examples. 

\section{ Failure of gap rigidity on the bidisk}

\subsection{Statement of main result and basic reductions}
\noindent For a product
domain $D^k$ denote by $\delta: D \to D^k$
the diagonal embedding given by $\delta(z) = (z, \cdots,z)$.
By Eyssidieux-Mok \cite{EM}, for any bounded
symmetric domain $D$, gap rigidity in the Zariski
topology holds for $(D^k,D;\delta)$.  It turns out
that the analogue is in general not valid when $\delta$
is replaced by another holomorphic totally geodesic embedding.

The main result of this section is to give a counterexample to the 
gap phenomenon as formulated in Eyssidieux-Mok \cite{EM} in this context.

\begin{theo}
Gap rigidity fails on the bidisk for $(\Delta^2,\Delta\times\{0\}))$.
\end{theo}

\noindent Since any bounded symmetric domain $\Omega$ of rank $\geq 2$ 
contains a totally
geodesic bidisk, Theorem 1 implies that gap rigidity fails in general whenever
$\Omega$ is of rank $\geq 2$.  We are going to prove Theorem 1 by constructing 
compact
holomorphic curves on products of compact Riemann surfaces of genus $\geq 2$
which are graphs of holomorphic maps.  The requirement that they become 
$\epsilon$-geodesic and modeled on $(\Delta^2,\Delta\times\{0\}))$ 
then translates to a uniform bound on the differential of the holomorphic 
map in terms of Poincar\'e metrics, as follows.

\begin{theo}
There exist sequences of compact Riemann surfaces $S_i$ and $T_i$ of 
genus $\geq 2$, together with a sequence of branched double covers 
$f_i: S_i \to T_i$, such that, writing $ds^2_{C}$ for the Poincar\'e 
metric of Gaussian curvature $-2$ on a compact Riemann surface $C$ of genus
$\geq 2$, and putting 
$\mu_i := 
\sup\Big\{\frac{\phantom{,}f^*_ids_{T_i}^2(x)\phantom{,}}{ds_{S_i}^2(x)}:
x \in S_i\Big\}$, we have $\lim_{i\to \infty} \mu_i = 0$.
\end{theo}

\noindent We now prove that Theorem 2 implies Theorem 1.

\begin{prv}
Let $G_i \subset S_i \times T_i$ be the graph of $f_i$.  Let
$F_i: \Delta \to \Delta$ be
an \text{arbitrary} lifting of $f_i: S_i \to T_i$ to the universal 
covering such that $F_i(0) = 0$.  Let 
$\widetilde{G}_i \subset \Delta^2$ be the graph of $F_i$, which 
is a lifting of $G_i \subset S_i \times T_i$ to the universal 
cover $\Delta^2$ of $S_i \times T_i$.  Since $F_i^*ds_{\Delta}^2
\leq \mu_i \cdot ds_{\Delta}^2$ and $\mu_i \to 0$ as $i \to \infty$, 
the holomorphic maps $F_i: \Delta \to \Delta$ converge to a
constant map $F_{\infty}$, $F_{\infty}(z) = 0$ for any $z \in \Delta$.
For $i \geq 1$, let $p_i \in \widetilde{G}_i$ be an arbitrary point.
Equip $\Delta^2$ with the product metric $ds_{\Delta^2}^2$
given by the Poincar\'e metric
$ds_{\Delta}^2$ in each factor, and denote by $B(p;r)$ the geodesic 
ball centred at $p$ and of radius $r > 0$ on $(\Delta^2;ds_{\Delta}^2)$. 
Consider the sequence of complex submanifolds $\widetilde{G}_i \cap B(p_i;1)
\subset B(p_i;1)$.  Let $\varphi_i \in \text{Aut}_o(\Delta^2)$ for which
$\varphi_i(p_i) = 0$.  Then, $\varphi_i(\widetilde{G}_i) \subset \Delta^2$
is again the graph of a holomorphic mapping $F'_i: \Delta \to \Delta$
which is a lifting of $f_i: S_i \to T_i$ satisfying $F'_i(0) = 0$.  The 
convergence argument applies therefore equally to the sequence $F'_i$
to show that  $\varphi_i(\widetilde{G}_i) \subset B(0;1)$ converge as 
complex-analytic subvarieties to $B_1 \times \{0\}$, 
where $B_1 \subset \Delta$ is the geodesic ball centred at $0$ of radius 1.
It follows that $G_i$ is $\epsilon_i$-pinched with $\epsilon_i$ converging 
to $0$, and $G_i \subset S_i \times T_i$ is modelled on 
$(\Delta^2;\Delta\times\{0\})$
for $i$ sufficiently large, in the sense of Eyssidieux-Mok \cite{EM}.  In
other words, gap rigidity fails for $(\Delta^2;\Delta\times\{0\})$.
\quad $\square{}$ \end{prv}

\noindent Let us give the idea of the proof of Theorem 2.
Let $f: S \to T$ be a
generically 
$s$-to-1 
holomorphic 
map 
between 
compact Riemann surfaces of genus $g(S), g(T) \geq 2$, and let $r$ be the
degree of the ramification divisor.   By the Riemann-Hurwitz formula we 
have $2(g(S) - 1) = 2s(g(T)-1) + r.$  Denote by $||\cdot||$ norms induced
by the Poncar\'e metrics $ds^2_S$ resp. $ds^2_T$ on $S$ resp. $T$ of 
constant curvature $-2$.  Then, $\int_S \|df\|^2 = s\cdot 
\text{Area}(T,ds^2_T)$.
It follows from the Gauss-Bonnet Theorem that the 
the average of $\|df\|^2$ over $S$ is $s\cdot 
\frac{\text{Area}(T,ds^2_T)}{\text{Area}(S,ds^2_S)}
= s\cdot\frac{g(T)-1}{g(S)-1}$, 
which is $\frac{2s(g(T)-1)}{2s(g(T)-1) + r}$ by the 
Riemann-Hurwitz formula.  If we fix $T$ and $s > 1$ then 
this average would get arbitrarily small by choosing $r$ arbitrarily 
large.  It is however not clear at all how sup$||df||$ can be bounded in 
terms of the ramification divisor.  We will fix some compact Riemann
surface $T$ and take double covers, and our approach is to choose the 
branching locus to be more or less evenly spaced on $T$.  
One can most properly talk about an ``evenly spaced'' set of points on 
an elliptic curve by exploiting the group structure.  For instance, 
for integers $n > 1$ the set of $n$-torsion points can be considered
such a set.  We will take a ramified 
double cover $h: T \to E$ to get $T$ of genus $\geq 2$ and obtain 
$f_i: S_i \to T$ by taking double covers over $T$ such that the branching
locus of $h\circ f_i: S_i \to E$ is evenly spaced in a precise way,
allowing us to work on the elliptic curve $E$ by descent.
In order to get an estimate on sup$||df_i||$ we will exploit the 
group structure on $E$ and the uniqueness properties of singular 
Hermitian metrics of negative Gaussian curvature on $E$ with prescribed
(fractional) orders of poles at a finite number of points.

\subsection{Proof of Theorem 2}  \noindent We have the following standard 
Lemma.

\begin{lem}
Let $C$ be a compact Riemann surface, $n$ be a positive integer, and 
$q_1, \cdots, q_{2n}$ be distinct points on $C$.  Then, there exists a
compact Riemann surface $K$, and a branched double covering map 
$f: K \to C$ which is branched precisely over $q_1, \cdots, q_{2n}$.
In other words, $f^{-1}(q_k) = \{p_k\}$ and $f$ is precisely ramified
(with ramification order 1 { \rm{a fortiori}}) at the $2n$ distinct points 
$p_1, \cdots, p_{2n}$. 
\end{lem}

\begin{prv}
Let $[D]$ be the divisor line bundle of degree $2n$ defined by the 
set $\{q_1,\ldots,q_{2n}\}$ of simple points, and let $L$ be a holomorphic
line bundle such that $L^2 \cong [D]$.  $L$ exists because any point in 
$Pic^0(C)$ is 2-divisible.  Let now $s$ be the 
canonical 
section of the divisor line bundle $[D]$ and $t$ be the 2-valent section 
of $\pi: L \to C$ over $C$ whose graph is defined by 
$K = \{v \in L: v^2 = s(\pi(v))\}$.  Then $K \subset L$ is nonsingular,
and, putting $f = \pi|_K$, we have a double cover $f: K \to C$ branched 
precisely over the $2n$ distinct points $q_1, \cdots, q_{2n}$. \quad 
$\square{}$
\end{prv}

\noindent We will now give a proof of Theorem 2 by constructing an example
where all the target Riemann surfaces
$T_i$ are 
identical. 
 
\noindent Let
$L 
\subset 
\Bbb 
C$ 
be 
any 
lattice 
and 
denote 
by $E = \Bbb C/L$ the quotient elliptic curve.  Let $e \in E$ be a nonzero
torsion point of order 2.  By Lemma 1 there is a compact Riemann surface
$T$ and a double cover $h: T \to E$ branched precisely over 0 and $e$.
Write $q_1 = h^{-1}(0)$ and $q_2 = h^{-1}(e)$.  Let $m = 2i - 1$ be an 
odd positive integer and consider the holomorphic map $\Phi_m: E \to E$ defined
by $\Phi_m(x) = mx$. Let $D_i := \Phi_m^{-1}(\{0,e\}) \subset E$.  Note
that $D_1 = \{0,e\}$ and $\text{Card}(D_i) = 2m^2$.  Since $(m-1)e = i(2e) = 
0$ 
on $E$, so that $m\cdot e = e$, we have also $0, e \in D_i$ for each 
positive integer $i$.  Let $i \geq 2$, i.e., $m \geq 3$.  Again by Lemma 1
there is a compact Riemann surface $S_i$ and a double cover 
$f_i: S_i \to T$ branched precisely over the reduced divisor 
$\Delta_i := h^{-1}(D_i - D_1)$, noting that $\text{Card}(\Delta_i) = 
4(m^2-1)$.   
We claim that $f_i: S_i \to T$ gives a sequence of holomorphic maps such 
that, writing $\mu_i = \sup\Big\{\frac{\phantom{,}f_i^* 
ds_{T}^2(x)\phantom{,}}{ds_{S_i}^2(x)}:
x \in S_i\Big\}$, we have $\lim_{i\to \infty} \mu_i = 0$.

\noindent We proceed to compare the Hermitian metrics $ds_{S_i}^2$ and
$f_i^*ds_T^2$ by descending to the elliptic curve $E$.
Since $h: T \to E$ is a double cover, the holomorphic mapping 
on $T - h^{-1}(\{0,e\})$ switching the two points of fibers of $h$
extends to an automorphism $\sigma$ of $T$ fixing the ramification 
points.  $\sigma$ fixes the Poincar\'e metric $ds_T^2$ on $T$ of 
Gaussian curvature $-2$.  Hence, $ds_T^2$ descends to a Hermitian metric 
$\theta$
of Gaussian curvature $-2$ on $E - \{0,e\}$.  Extending across the 
branching points 0 and $e$, $\theta$ can be interpreted as a Hermitian
metric on $E$ with simple poles at $0$ and $e$.  More precisely, in terms
of a local holomorphic coordinate $z$ at 
$0$ 
or at 
$e$, we 
can 
write
$\theta = \frac{\phantom{,}a(z)}{|z|}|dz|^2$, where $a(z)$ is a continuous 
positive 
function which is smooth except at $z=0$.  We may say that $\theta$ is a 
Hermitian metric on 
the $\Bbb Q$-line bundle $T_E \otimes [D_1]^{-\frac{1}{2}}$.       
Similarly for $i > 1$, by means of the double cover 
$f_i: S_i \to T$, the Poincar\'e metrics $ds_{S_i}^2$ of Gaussian 
curvature $-2$ descends to a Hermitian metric $\eta_i$ on the 
$\Bbb Q$-line bundle $T_T \otimes [\Delta_i]^{-\frac{1}{2}}$.  We 
observe that $\eta_i$ is invariant under the holomorphic involution
$\sigma$ on $T$.  In fact $f_i^*(\sigma^*\eta_i)$ gives a smooth 
Hermitian metric on $S_i$ of constant Gaussian curvature $-2$, so
that $f_i^*(\sigma^*\eta_i) = ds_{S_i}^2$ by the uniqueness of 
Poincar\'e metrics of Gaussian curvature $-2$, i.e., by the Ahlfors-Schwarz
Lemma.  It follows that $\sigma^*\eta_i = \eta_i$, as observed.
As a consequence $\theta_i$, $i > 1$, descends to a Hermitian metric 
$\theta_i$ on 
$E - D_i$.  Across points of $D_i - D_1, \theta_i$ extends to a Hermitian
metric with a simple pole since $h: T \to E$ is unramified over 
$\Delta_i = h^{-1}(D_i-D_1)$.  On the other hand, $\theta_i$ extends
across points of $D_1 = \{0,e\}$ as a Hermitian metric with a simple
pole because $\eta_i$ is smooth and positive at the two points
$h^{-1}(0) = q_1$ and $h^{-1}(e) = q_2$, and $h$ is ramified at
$q_1$ resp. $q_2$ to the order 1. More precisely, $\theta_i$ 
is a Hermitian metric on the $\Bbb Q$-line bundle
$T_E \times [D_i]^{-\frac{1}{2}}$ for $i>1$.  
We will also write $\theta_1$ for $\theta$, 
so that the last statement is true for $i \geq 1$.

\noindent For $i > 1$ and for any $x \in S_i$,
$\frac{\phantom{,}f_i^*ds_{T}^2(x)\phantom{,}}{ds_{S_i}^2(x)}$ 
is the same as $\frac{\phantom{,}\theta_1(y)\phantom{,}}{\theta_i(y)}$, 
where $y = h\circ f_i(x)$, 
provided that $y \notin D_i$.  For such points $y \in E$ the 
inequality  $f_i^*ds_{T}^2(x) \leq ds_{S_i}^2(x)$ translates into 
the inequality $\theta_1(y) \leq \theta_i(y)$.  
Since the endomorphism $\Phi_m:E\to E$ is unramified, $\Phi_m^*(\theta_1)$
gives a smooth Hermitian metric on the $\Bbb Q$-line bundle 
$T_E \otimes [D_i]^{-\frac{1}{2}}$.  As in the last paragraph, by 
pulling-back to the double cover $R_i$ over $E$ branched precisely over
$D_m$, $\text{Card}(D_m) = 2m^2$, from the uniqueness of the Poincar\'e
metric of Gauss curvature $-2$ on $R_i$ it follows that $\theta_i =
\Phi_m^*\theta_1$.  We have therefore
$\mu_i = 
\sup\Big\{\frac{\phantom{,}\theta_1(y)\phantom{,}}{\Phi_m^*\theta_1(y)}: 
y \notin D_i\Big\}$, where $m = 2i -1$.
It remains to show that $\mu_i \to 0$ as $i \to \infty$.

\noindent Fix a Euclidean metric $\omega$ on $E$ such that $\theta_1 > \omega$.
Then $\Phi_m^*\theta_1 > m^2\cdot \omega$.  Fix a coordinate unit disk $U_1$
at $0$ (resp. $U_2$ at $e$), with coordinte $z$, so that $z(0) = 0$
(resp. $z(e) = 0)$, and such that $\theta_1 = \frac{a(z)}{|z|}|dz|^2$
with $a(z)$ a continuous positive function bounded between two positive 
constants.  There is a positive constant $K$ such that $\theta_1(y)
\leq K\omega$ for $y \in E - U_1 - U_2$.  Hence 
$\frac{\theta_1(y)}{\phantom{,}\Phi_m^*\theta_1(y)\phantom{,}} < 
\frac{K}{\phantom{,}m^2\phantom{,}}$ 
for $y \in E-U_1-U_2$. At $0 \in E$ we may choose $U_1$ such that $\Phi_m(z) = 
m\cdot z$
as a germ at 0. For $|z| < \frac{\phantom{,}1\phantom{,}}{m}$,
$\Phi_m^*\theta_1(z) = \frac{a(mz)}{\phantom{,}|m\cdot z|\phantom{,}}m^2|dz|^2
= \frac{\phantom{,}m\cdot a(mz)\phantom{,}}{|z|}|dz|^2$, so that 
$\frac{\theta_1(z)}{\phantom{,}\Phi_m^*\theta_1(z)\phantom{,}}
< \frac{a(z)}{\phantom{,}m\cdot a(mz)\phantom{,}}
< \frac{\phantom{,}C_1\phantom{,}}{m}$ for some constant $C_1$.  
On the other hand, for $|z| > \frac{1}{\phantom{,}m\phantom{,}}$ 
we have $\theta_1(z) < C_2m\cdot \omega$ for some positive
constant $C_2$, while $\Phi_m^*\theta_1(z) > m^2\cdot \omega$, so 
that $\frac{\theta_1(z)}{\Phi_m^*\theta_1(z)} < 
\frac{\phantom{,}C_2\phantom{,}}{m}$.  
We have therefore
$\frac{\theta_1(z)}{\phantom{,}\Phi_m^*\theta_1(z)\phantom{,}} 
< \frac{\phantom{,}C\phantom{,}}{m}$ 
for $C = \text{max}(C_1,C_2)$,
for $z \in U_1$ and similarly for $z \in U_2$.  Combining with the estimate
on $E-U_1-U_2$, we have for $m = 2i-1$ sufficiently large the estimate
$\frac{\phantom{,}\theta_1(y)\phantom{,}}{\Phi_m^*\theta_1(y)} 
< \frac{\phantom{,}C\phantom{,}}{m}$ for all $y \in E-D_i$.
It follows that $\mu_i < \frac{\phantom{,}C\phantom{,}}{m} \to 0$ as $i \to 
\infty$.
The proof of Theorems 1 and 2 is complete. $\quad\quad \square$. 

{\bf Remarks.}  We note that the unit disk plays a very
special role in the proof of Theorem 1.  The analogue
of Theorem 1 does not apply when $\Delta$ is replaced
by an irreducible bounded symmetric domain $D$ of
rank $\geq 2$ (cf. (3.3), Proposition 4).  The case
where $D$ is of rank 1 but of dimension $\geq$ 2
remains open.

\noindent
\section{ Gap rigidity in the Zariski topology by intersection theory}

\noindent
\subsection{Basic facts and notations on geodesic embeddings of bounded 
symmetric domains}

\noindent
\subsubsection{Lie theoretic data attached to a symmetric domain}

\noindent There is a 1-1 correspondence between bounded symmetric domain
$(\Omega,o)$
and semisimple Lie algebras of Hermitian type. These are
 Lie theoretic data $({\mathfrak g},H_0)$, $(\mathfrak g,\theta)$ a semisimple 
Lie algebra
with a Cartan involution and
$H_0$ an element of the center of the associated maximal compact subalgebra 
such  that $\mathop{ad} (H_0)^2=\theta$ (see \cite{Sat}).

\noindent If $\Omega$ is irreducible,
the correspondence can be described as follows: $\Omega= G/K$ and $o=eK$ where
 $G$ is the Lie group underying the real points
of a connected almost simple real algebraic group, also denoted by $G$,
such that $G^{{\mathbb C}}=G\otimes _{{\mathbb R}} {\mathbb C}$
is connected and simply connected,
 $K$ is a maximal compact subgroup
such that the center of $K$ is isomorphic to $U(1)$.
The adjoint group of
$G$ is the identity component of $\mathop{Aut} (\Omega)$.
Let $\mathfrak g$ be the Lie algebra of $G$, $\mathfrak l$ be the Lie algebra 
of $K$, $\theta$ the Cartan involution of the symetric pair
$(G,K)$
, $\mathfrak g=\mathfrak l\oplus\mathfrak p$ the Cartan decomposition.
We choose  $H_0 \in {\mathfrak z}$ satisfying
$\mathop{ad} (H_0)^2=\theta$, such that, when restricted to $\mathfrak p$,
$\mathop{ad} (H_0)$ corresponds  to the almost complex structure operator of 
$\Omega$ under
the canonical isomorphism ${\rm can}: T^{{\mathbb R}}_o\Omega \to{\mathfrak 
p}$ . In 
particular
${\mathfrak p}^{{\mathbb C}}={\mathfrak p}^+\oplus{\mathfrak p}^-$, where
$\mathfrak p^{\pm}$ is the $\pm \sqrt{-1}$-eigenspace of $\mathop{ad}(H_0)$
and $\mathfrak{p}^+$ corresponds to $T^{1,0}\Omega$.
Let $(\ . \ ,\ . \ )_{{\mathcal K}}$ be the Killing form on $\mathfrak g$.
A $G$-invariant Riemannian metric $g^{\Omega}_o$ can 
be constructed such that
$g^{\Omega}_o|_{T_o\Omega}={\rm can}^* ( \ . \ ,\ . \ )_{{\mathcal K}}| 
_{{\mathfrak 
p}}$, $g^{\Omega}_o$ is 
K\"ahler-Einstein.

\noindent A general bounded symetric domain $(\Omega,o)$ splits as a product
of irreducible ones $\Omega=\Omega_1\times\ldots\times\Omega_a$
and we define ${\mathfrak g}={\mathfrak g}_1\times\ldots\times {\mathfrak 
g}_a$, $H_0= (H_0)_1\times \ldots \times (H_0)_a$,
$G_0= (G_0)_1\times \ldots \times (G_0)_a$,
$K=K_1\times\ldots\times  K_a$. 

\noindent The converse correspondence constructing from the data $({\mathfrak 
g}, H_0)$
as above
a bounded symmetric domain in a complex vector space is a celebrated theorem 
of 
Harish-Chandra's, the classical cases being
in E. Cartan's thesis.

\noindent The action of $K\subset G$ induced by the adjoint action of $G$ on 
$\mathfrak
g$ leaves invariant the decomposition ${\mathfrak g}^{{\mathbb C}}={\mathfrak 
p}^+\oplus{\mathfrak p}^-\oplus {\mathfrak l}^{{\mathbb C}} $ and gives rise
to the isotropy representation of $K$
on $T_o^{1,0}\Omega ={\mathfrak p}^+$.

\subsubsection{Isotropy action on the Grassmannian of p-planes in 
${\mathfrak p}^+$}

\noindent Consider ${\rm Gr}(p, T_o\Omega)$ the Grassmannian parametrizing 
complex
$p$-planes
in ${\mathfrak p}^+=T^{1,0}_o\Omega$.

\noindent We endow $\Lambda^p {\mathfrak p}^+$
with the Hermitian metric functorially attached to the Hermitian
metric on ${\mathfrak p}^+$ defined by the formula $(\alpha,\beta)_{{\mathfrak 
p}^+}=(\alpha,\overline{\beta})_{{\mathcal K}}$.
This defines a Fubini-Study K\"ahler metric on ${\mathbb P}(\Lambda^p 
{\mathfrak p}^+)$. Consider the K\"ahler
metric $ds^2_p$ on ${\rm Gr}(p, T_o\Omega)$ induced by this Fubini-Study 
metric 
under 
the Pl\"ucker embedding ${\rm Gr}(p, T_o\Omega) \to {\mathbb P}(\Lambda^p 
{\mathfrak 
p}^+)$. Call this the Fubini-Study metric on ${\rm Gr}(p, T_o\Omega)$. The 
$K$-invariant 
K\"ahler  metrics on ${\rm Gr}(p, T_o\Omega)$ take
 the form $const. ds^2_p$, if $\Omega$ is irreducible.

\noindent Let ${\mathcal A}\in {\rm Gr}(p, T_o\Omega)$ a $p$-plane in 
${\mathfrak 
p}^+$.
Let $B=(e_1, \ldots, e_p)$ a unitary basis of ${\mathcal A}$.
Let $S(B,{\mathcal A})= \sqrt{-1} \sum_{i=1}^p [e_i,\overline{e_i}]$, where 
$[-,-]$
is the Lie bracket in $\mathfrak g^{\mathbb C}$. This expression is obviously 
independant of $B$ and defines a real-analytic mapping $\Sigma :{\rm 
Gr}(p, T_o\Omega)
\to {\mathfrak l}$.

\subsubsection{Embeddings of symmetric domains}

  \noindent Consider  two  Lie algebras of Hermitian type $({\mathfrak g}, 
H_0)$ and
$({\mathfrak g}', H'_0)$. A $(H_1)$-homomorphism $\rho:({\mathfrak g}, H_0)
\to({\mathfrak g}', H'_0)$ is a Lie algebra morphism $\rho:{\mathfrak 
g}'\to{\mathfrak g}$ such that
$\mathop{ad}(H_0) \rho = \rho \mathop{ad}(H'_0)$. An $(H_2)$-homomorphism
is a Lie algebra morphism satisfying $\rho(H'_0)=H_0$.
$(H_2)$-homomorphisms are $(H_1)$, see  \cite{Sat}, pp. 83-88. 
These notions
have been introduced by Satake \cite{Sat1}.

\noindent Totally geodesic embeddings of pointed symmetric domains are in
1-1 correspondence with
injective
$(H_1)$-homomorphisms. Indeed,  such  a $\rho$ yields a morphism of the
underlying symmetric Lie algebras 
that is complex linear for the complex structures on ${\mathfrak p}'$
(resp. ${\mathfrak p}$) induced by ${\mathop{ad} } (H_0)$ and can be 
exponentiated to an
homomorphism $\rho: G'\to G$ satisfying $\rho(K')=K$, so that the 
map $G'/K' \to G/K$ is holomorphic.
Those corresponding to $(H_2)$-homomorphisms will be called $(H_2)$.

\begin{lem}\label{sig}
Let $\Omega$ be irreducible. Then there exists a real number 
$c_{\Omega}$ such that $\Sigma({\mathfrak p}^+)= \sqrt{-1} c_{\Omega}
 H_0$.
\end{lem}

\begin{prv}
The lemma follows directly from independance of $\Sigma$ with the respect to 
the choice of a basis. Indeed, we may  use a change of basis which lies in  
$K$ to see that $\Sigma$ is a fixed point of the adjoint action
of $K$. Hence the lemma.

\end{prv}

\noindent
Let $\Omega'$ be an irreducible subdomain of $\Omega$ and
$\rho: {\mathfrak g}'\to {\mathfrak g}$ be the corresponding $(H_1)$-embedding.
Let $d_{\Omega',\Omega}>0$ be the real number defined by
$g^{\Omega}_o|_{\Omega'}= d_{\Omega',\Omega}. g_o^{\Omega '}$.

\noindent 
In order to  state our results on gap rigidity in the Zariski topology, we 
need to introduce a variant of the 
$(H_1)$ and $(H_2)$ conditions.
Consider an embedding of
bounded symmetric domains $j:\Omega'\to \Omega$ and
$\rho: {\mathfrak g}'\to {\mathfrak g}$ be the corresponding $(H_1)$-embedding.
Let $\Omega'=\Omega_1'\times \ldots \times \Omega'_a$
be the irreducible decomposition of $\Omega'$. We will say 
$j$, resp. $\rho$, is $(H_3)$ iff the following holds:
$$
\rho (\sum_{i=1}^a c_{\Omega'_i} d_{\Omega'_i,\Omega} H'_{0i}) \in {\mathbb R}
H_0.
$$
 \begin{lem} $(H_3)$-embeddings are $(H_2)$. A $(H_2)$-embedding
  $\Omega'\to \Omega$ is
 $(H_3)$ iff the Einstein constants of the metrics $g_o^{\Omega}|_{\Omega'_i}$
are the same, where $\Omega'=\Omega'_1\times \ldots \times \Omega'_a$ 
and $\Omega'_i$ is irreducible.
 \end{lem}
\begin{prv}
Let $H'_1=A.\sum_{i=1}^a c_{\Omega_i'} d_{\Omega_i',\Omega} H'_{0i}$, $A$ a 
constant
such that $\rho(H'_1)=H_0$. For every $z\in{\mathfrak g}'$, we have since
$\rho$ is $(H_1)$, $[H_0-\rho(H'_0), \rho (z)]=0$,
 hence $[\rho(H'_1-H'_0),\rho(z)]=0$ and $[H'_1-H'_0,z]=0$. 
 Since ${\mathfrak g}'$ is semisimple this implies that $H_1'=H_0'$.
 A fortiori, $\rho$ is $(H_2)$. It also follows that for every $1\le i \le 
a$,
 $A.c_{\Omega_i'} d_{\Omega_i',\Omega}=1$. A standard curvature formula implies
that
 $c_{\Omega_i'} d_{\Omega_i',\Omega}$ is minus 
 the Einstein constant of $g_o^{\Omega}|_{\Omega'_i}$.
\end{prv}

\noindent In particular, for an embedding of irreducible
bounded symmetric domains $\Omega'\subset \Omega$, the $(H_2)$ and $(H_3)$
conditions are equivalent.

\subsection{Statement of main results}
\noindent The purpose of this section is to give
a general criterion under 
which gap rigidity holds for a pair $(\Omega,D)$ in the Zariski topology
with $\Omega$ irreducible (and of rank $\geq 2$). 

\noindent Let us first state a GIT interpretation of the $(H_3)$-condition. Say
${\mathcal A}\in {\rm Gr}(p, T_o\Omega)$ is GIT-semistable if there is a
$K$-invariant 
complex closed hypersurface $\mathcal Z_o \subset {\rm Gr}(p, T_o\Omega)$ such 
that 
${\mathcal A}\notin \mathcal Z_o$.

  \begin{prop} \label{sth2}
  The embeddings of symmetric domains
   $(\Omega',o)\to (\Omega,o)$  such that
   $T_o\Omega'$ is GIT-semistable
   in ${\rm Gr}({\dim\Omega'}, T_o\Omega)$ with respect to the isotropy action 
of 
$K$
   are precisely the $(H_3)$-embeddings.
  \end{prop}

\noindent
As a step towards his classification  of embeddings of symmetric domains,
  Ihara \cite{Iha} obtained the full classification theory of
$(H_2)$-embeddings.

\noindent 
It turns out that, thanks to Proposition \ref{sth2}, images of 
$(H_3)$-embeddings of
bounded symmetric domains play an important role
in the question of gap rigidity in the Zariski
topology.  A compact totally-geodesic
complex submanifold $S$ of a quotient $X$ of a bounded
symmetric domain by a discrete group of biholomorphic
automorphisms will be referred to as a holomorphic
geodesic cycle.  If $S \subset X$ arises from an
$(H_3)$-embedding, $S$ will be referred to as an
$(H_3)$-holomorphic geodesic cycle.  We prove:

\begin{theo}
Let $\Omega = G/K$ be an $n$-dimensional irreducible  bounded symmetric
domain.  Let 
$\Gamma$ be a torsion-free discrete group of biholomorphic 
automorphisms of $\Omega$, and write $X := \Omega/\Gamma.$
Assume there exists $D \subset \Omega$ an $(H_3)$-embedding, $o \in D$,
$\dim(D):= p$. 

 Fix a projective $K$-invariant
hypersurface $\mathcal 
Z_o 
\subset 
\text{\rm 
Gr}(p, T_o\Omega)$ 
such 
that 
$[T_o(D)]
\notin \mathcal Z_o$.  Denote by $\pi: \Bbb 
G(X) 
\to 
X$ 
the 
Grassmann 
bundle 
of $p$-dimensional tangent planes and write $\mathcal 
Z = 
\mathcal 
Z_X 
\subset 
\Bbb 
G(X)$  
for the locally homogeneous 
subbundle 
of 
projective 
hypersurfaces
corresponding to $\mathcal
Z_o $. 

Let $S \subset X$ be a compact complex p-dimensionnal submanifold such that 
for any $x \in S$,
$[T_x(S)] \notin \mathcal 
Z_x$.  
Then, 
$S 
\subset 
X$ 
is an $(H_3)$-holomorphic
geodesic 
cycle.
\end{theo}

\noindent {\bf Remark}  We do not exclude the possibility that $S$ is 
uniformized
by some $(H_3)$ complex totally-geodesic submanifold $D' \subset \Omega$
which is not equivalent to $D$ under Aut$(\Omega)$.  This may in fact happen,
but
there are up to equivalence under Aut$(\Omega)$ only a finite number of
possibilities.

\subsection{Proof of Proposition \ref{sth2}}

\subsubsection{Moment map of the isotropy action on the Grassmannian of planes 
in ${\mathfrak p}^+$}

  \noindent
Let $\omega_p$ be the K\"ahler form associated to $ds^2_p$. View it as a 
symplectic form.
$K$ acts on the symplectic manifold
 $({\rm Gr}(p, T_o\Omega), \omega_p)$ preserving the symplectic form. Since 
the 
Grassmannian is simply connected, there exists a moment map $\mu_p$
 for this symplectic action \cite{Sou} (see \cite{GIT3}, Chap. 8).

  \noindent
 It might be useful to recall the definition of this central concept of 
symplectic geometry.
 Let $(S,\omega)$ be symplectic manifold acted upon by a
 connected  Lie group $M$ whose Lie algebra is denoted by $\mathfrak m$, a   
map $\mu: S\to {\mathfrak m}^*$
 is called a moment map if it is smooth, $M$-equivariant with respect to
 the given action of $M$ on $X$ and the coadjoint action and

 $$ \forall x\in M, \xi \in T_xS, l\in {\mathfrak m}, \ \ << d\mu(x), \xi>,
l>= <\omega(x), \xi\wedge l_x>
 $$

 \noindent where $<-,->$ denotes
 the canonical pairing between a vector space $E$ and its dual $E^*$ and 
$l_x=\frac{d}{dt} exp(tl)x  |_{t=0}$.

\noindent When it exists, a moment map is unique up to the addition of a 
coadjoint fixed point in
${\mathfrak m}^*$. When $M={\mathbb R}$, a moment map is a global hamiltonian
function for the flow of $\frac d {dt}$.

\noindent If $V$ is a complex vector space with a non degenerate sesquilinear
pairing $h$, we also denote by $h$ the corresponding conjugate linear
isomorphism $V\to V^*$. In the next lemma, $V={\mathfrak l}^{\mathbb C}$
and $h_{{\mathfrak g},{\mathfrak l}}(\lambda, \mu)= ( \lambda, 
\overline{\mu})_{{\mathcal K}}$.

\begin{lem}\label{moment}
The map $\mu_p= h_{{\mathfrak g},{\mathfrak l}} \circ \Sigma$ is a moment map 
$\mu_p: {\rm Gr}(p, T_o\Omega) \to {\mathfrak l}^*$ relative to the symplectic 
action of 
$K$
on $({\rm Gr}(p, T_o\Omega),\omega_p)$.
\end{lem}

\subsubsection{Proof of Lemma \ref{moment}}

\noindent

\begin{lem}
Let the Lie group $L_1\times L_2$ act symplectically on $(S,\omega)$
with moment map $\mu=(\mu_1, \mu_2): S\to {\mathfrak l}^*_1
\times {\mathfrak l}^*_2$. Let $z$ be a coadjoint fixed point of $L_2$
which is furthermore a regular value of $\mu_2$ and assume $L_1$ acts freely 
on $\mu_2^{-1}(z)$. The symplectic
quotient \cite{MW} $S//L_2=\mu_2^{-1}(z)/ L_2$ is acted upon symplectically by 
$L_1$.
Let $\pi: \mu_2^{-1}(z) \to S//L_2$ be the canonical quotient map.

 The symplectic action of $L_1$ on $S//L_2$ admits a moment map $\mu$
 defined by the relation $\mu\circ \pi=\mu_1| _{\mu_2^{-1}(z)}$.
\end{lem}

\begin{prv} $L_1$-equivariance is clear.
Recall (see \cite{GIT3}, p. 146) that $\mu_2^{-1}(z)$ is a coisotropic 
submanifold of $S$ whose isotropic foliation $\ker (\omega|_{\mu_2^{-1}(z)})$ 
is precisely the foliation by $L_2$-orbits
and that the
symplectic form on $S//L_2$ is the form induced by $\omega$ on local 
transversal sections of $\pi$. The required differential identity follows then 
directly from the definitions.
\end{prv}

\noindent An easy example is the moment map $\mu_n$ associated with the action
of $U(n)$ on ${\mathbb C}^n$ with coordinates $(z^i)$ equipped with 
$\frac{\sqrt{-1} }{2\pi}
\sum_i dz^i\wedge d\bar z^i$. Identifying $\mathfrak{u}(n)$ and its dual
by means of the scalar product $(x,y)_{\mathfrak{u}(n)}=-tr(xy)$. We have:

 $$\mu_n((z^1, \ldots, z^n) )= -\frac{ \sqrt{-1}}{ 2\pi}(z^i \bar z^j)_{1\le 
i,j\le n}.$$

 \noindent Since the moment map for the action of $K$ on a product $X\times Y$
 of two symplectic $K$-manifolds is the sum of the moment maps for $X$ and 
$Y$, it follows that the moment map for the action of $U(n)\times U(p)$
 on complex matrices $Z=(z^i_l)_{1\le i \le n, 1\le l \le p}$
 with $n$ rows and $p$ columns equipped with
 the Hermitian scalar product $Z\mapsto\frac{1}{\pi} Tr (Z ^t\bar Z)$ is 
$\mu(Z) = -\frac{ \sqrt{-1}}{ 2\pi}(Z ^t\bar Z, -^t\bar Z Z)=-\frac{ 
\sqrt{-1}}{ 2\pi}((\sum_{l=1}^p z^i_l \bar z^j_p),-(\sum_{i=1}^n z^i_l \bar 
z^i_m)).$

 \noindent If $p\le n$ the Grassmannian of $p$ planes in ${\mathbb C}^n$
equipped with its Fubini-Study 2-form is the
 symplectic quotient $M_{n,p}$ by $U(p)$
corresponding to the regular value $\frac{ \sqrt{-1}}{ 2\pi} Id_p$.

\noindent This gives the expression $\mu([Z])= -\frac{ \sqrt{-1}}{ 2\pi} Z 
^t\bar Z$
for any matrix $Z$ representing a unitary basis of $[Z]$ in this Grassmannian.

\noindent Let us restate this fact in slightly more invariant terms.
Let $(V,h)$ be a Hermitian vector space.
Use the canonical complex linear identifications
${\mathfrak u}(V,h)^{\mathbb C}=\mathop{End}(V) \simeq V\otimes V^*$,
${\mathfrak u}^*(V,h)^{\mathbb C}\simeq V^*\otimes V$ and
the conjugate linear isomorphism $h:V\to V^*$.
Let $\mathcal A$ a $p$-plane in $V$ with unitary basis
$(v_1,\ldots,v_p)$. The value at $\mathcal A$  of the moment map with respect 
to the Fubini Study metric is $\mu_p({\mathcal A})=
\frac {\sqrt{-1}}{2\pi} \sum_{i=1}^p h(v_i) \otimes v_i$.

\noindent Obviously, if $(S,\omega)$ is symplectic manifold acted upon by a
 connected  Lie group $M$  admitting a moment map and $M'\subset M$ is a Lie 
subgroup the moment map for $M'$ is the composition of the moment map for $M$
 with the canonical map $ \mathfrak m^* \to ({\mathfrak m}')^*$.

\noindent We now compute the map $ {\mathfrak p}^{+*}\otimes {\mathfrak p}^{+}
\simeq {\mathfrak u}^*((-,-)_{{\mathfrak p}^+})^{\mathbb C} \to {\mathfrak
l}^{{\mathbb C }*}$ induced by the isotropy representation $\rho:K\to 
U((-,-)_{{\mathfrak
p}^+}) $. Fix $v\in {\mathfrak p}^+$, $l\in {\mathfrak l}$ and $(e_i)$ a 
unitary basis of this vector space. We have:

\begin{eqnarray*}
< h_{{\mathfrak p}^+}(v)\otimes v, \rho(l)>&=&
<h_{{\mathfrak p}^+}(v)\otimes v, \sum_i [l,e_i] \otimes h_{{\mathfrak 
p}^+}(e_i)>\\
&=&
\sum_i < h_{{\mathfrak p}^+}(v), [l,e_i]><v, h_{{\mathfrak p}^+}(e_i)>\\
&=&
<h_{{\mathfrak p}^+}(v),[l,v]>= ([l,v],v)_{{\mathfrak p}^+}
=([l,v],\bar v)_{{\mathcal K}}\\
&=&(l,[v,\bar v])_{{\mathcal K}} .
\end{eqnarray*}

\noindent The last equality holds since the adjoint representation acts by  
isometries
of the Killing form.
This concludes the proof of Lemma \ref{moment}.

\subsubsection{GIT analysis of embeddings of bounded symmetric domains}

\paragraph{Semistability of $(H_3)$-embeddings}

\noindent Lemma \ref{moment},  Lemma \ref{sig} above and Kempf-Ness theory,
 a classical relation between moment maps and Geometric invariant theory
 (Theorem 8.3 of \cite{GIT3}), imply

\begin{coro} \label{corol}
Let ${\mathcal A} \in {\rm Gr}(p, T_o\Omega)$ such that $\Sigma(\mathcal A)= c
H_0$.
Then ${\mathcal A}$ is GIT semistable with respect to the action of 
$K^{\mathbb C}$.
\end{coro}

\noindent The $K$-action reduces to a $K/U(1)$ on the Grassmannian because
$U(1)$ acts by complex homotheties on
${\mathfrak p }^+$. In particular the moment map composed with
the projection to ${\mathfrak u}(1)^*$ is constant.
So that the condition $\Sigma(\mathcal A)= c H_0$ means the vanishing of
the moment map for $K/U(1)$ on ${\rm Gr}(p, T_o\Omega)$.

\noindent Alternatively,
the values of the moment map are in a fixed affine hyperplane $H$ of the form
$H=\{ \zeta \in
{\mathfrak l}^*, \ {\rm such \ that} \ <\zeta,H_0> ={\rm cste}\}$ . Using
lemma \ref{moment},  we see that for an aribitrary subspace ${\mathcal B}$,
$\Sigma({\mathcal B})=c H_0 + \Sigma'({\mathcal B})$  with
$h_{{\mathfrak g},{\mathfrak l}}( H_0, \Sigma'({\mathcal B}))=0$.
For future use, we record the following observation

\begin{lem} \label{absmin}
Let ${\mathcal A} \in {\rm Gr}(p, T_o\Omega)$ such that $\Sigma(\mathcal A)= c
H_0$.   $\| \Sigma \|^2=\|
\Sigma\|_{{\mathfrak g},{\mathfrak l}}^2$ achieves its absolute minimum
 at ${\mathcal A}$.
\end{lem}

\noindent
Corollary \ref{corol}  gives one implication in  Proposition \ref{sth2}.
We now prove the converse
statement.

\paragraph{Embeddings of symmetric domains as critical points of the
\lq Morse function\rq \ attached to the moment map}

\begin{lem}\label{crit}
Let $(\Omega',o)\to (\Omega,o)$ be a geodesic embedding. Then
$[\rho({\mathfrak p}^+)]$ is a critical point of $\| \Sigma \|^2$.
\end{lem}
\begin{prv}  The tangent space at ${\mathcal A}\in {\rm Gr}(p, T_o\Omega)$
is isometric to $H=Hom({\mathcal A},{\mathcal A}^{\perp})$.
If $v\in H$ and $(e_i)$ is a unitary basis of ${\mathcal A}$, $(e_i+t.v(e_i))$
is unitary up to the second order term in $t$ thus
$\partial_v \|\Sigma ({\mathcal A})\|^2= \sum_{ij} ([v(e_i),\bar e_i], 
[e_j,\bar e_j])+ (i\leftrightarrow j)= 2 \sum_i (v(e_i), [\Sigma,e_i])$. Hence 
${\mathcal A}$ is critical iff
$[\Sigma({\mathcal A},{\mathcal A})]\in {\mathcal A}$.
When ${\mathcal A}=\rho({\mathfrak p'}^+)$ we have actually
$\forall \alpha,\beta, \gamma \in {\mathcal A}, [[\alpha,\bar \beta],\gamma] 
\in {\mathcal A}$.
\end{prv}

\paragraph{Conclusion of the proof}
\noindent Kempf-Ness theory can be made more precise than Theorem 8.3 of 
\cite{GIT3}.
The statement we need is thm 8.10 p. 109 in \cite{Kir} through the following
consequence:

\begin{prop}. Let $K$ be a compact Lie group.
Assume $K^{\mathbb C}$ acts linearily on a projective manifold $X$ embedded
in ${\mathbb P}^N$.
Let $\mu_{{\mathbb P}U(N+1)}: X \to {\mathfrak{su}}(N+1)^{*}$ be the moment 
map corresponding to a
the Fubini-Study metric and $\mu:X\to {\mathfrak l}^{*}$ be the moment map 
obtained by composing
$\mu_{{\mathbb P}U(N+1)}$ and the canonical surjection ${\mathfrak{su}}(N+1)^* 
\to {\mathfrak l}^{*}$.

The set of semistable points of $X$ is the minimal Morse stratum, that is the
set
of points attracted by $\mu^{-1}(0)$ under the steepest descent flow of 
$\|\mu\|^2$.

In particular, every  critical point of $\|\mu\|^2$ whose critical value is not
$0$ is unstable
in the GIT sense.
\end{prop}

\noindent Since the moment map we use has the required form,  applying Lemma 
\ref{crit}, we deduce Proposition \ref{sth2}.

\subsection{Proof of Theorem 3}

\noindent Denote also by $\pi: \Bbb G(\Omega) \to \Omega$ the Grassmann bundle 
of
$p$-planes in $T(\Omega)$, and by 
$\mathcal
Z_{\Omega} 
\subset 
\Bbb 
G(\Omega)$ 
the 
homogeneous bundle of projective hypersurfaces corresponding to 
$\mathcal Z_o \subset \text{\rm Gr}(p,T_o(\Omega)) = \Bbb 
G_o(\Omega)$.  
Let 
$\Omega \subset M$ be the 
Borel embedding of $\Omega$ into its compact dual $M = G_c/K$, 
and $N\supset D$ be the compact dual of $D$. $N$ is a complex
submanifold of $M$. The complexification
$K^{\Bbb C}$ of $K$ acts naturally as a group of automorphisms on $M$.  Since 
$\mathcal Z_o \subset \Bbb G_o(\Omega)$ is complex-analytic, it is 
invariant 
under
$K^{\Bbb C}$ and the image of 
$\mathcal
Z_o 
\subset 
\Bbb 
G_o(\Omega)$ 
under 
the 
action 
of Aut$(M) = G^{\Bbb C}$ defines a holomorphic subbundle 
$\mathcal 
Z_M 
\subset 
\Bbb G(M)$ of the Grassmann bundle of $p$-planes in $M$. 
$\mathcal 
Z_{\Omega}$ 
is 
the restriction of ${\mathcal
Z}_M$ 
to 
$\Omega$, 
and 
it 
descends 
to 
$\mathcal Z 
= 
\mathcal Z_X \subset \Bbb G(X)$ under the action of $\Gamma$.   The 
hypersurface
$\mathcal Z_M \subset \Bbb G(M)$ defines a $K$-invariant divisor line 
bundle.  
$\Bbb G(M)$ embeds 
canonically 
into 
$\Bbb
P(\Lambda^pT_M)$, 
where 
at 
each base point $x \in M$ the embedding $\Bbb G_o(M) \subset \Bbb 
P(\Lambda^pT_o(M))$
is congruent to the Pl\"ucker embedding.  The 
tautological 
line 
bundle 
$L$ on $\Bbb P(\Lambda^pT_M)$ restricts to the tautological 
line 
bundle on 
$\Bbb G(M)$, to be denoted also by $L$.
Since the Picard group of both $M$ and the typical 
fiber 
$\Bbb 
G_o$ of 
$\pi: \Bbb G(M) \to M$ are both infinite cyclic, the Picard group of the 
total space $\Bbb G(M)$ is isomorphic to $\Bbb Z^2$, and is generated as 
a group by $\pi^*{\mathcal 
O}(1)$ 
and by 
$L$, 
where 
${\mathcal 
O}(1)$ 
denotes 
the 
positive generator of $Pic(M)$.  Thus there exist positive integers $m, \ell$
such that $\mathcal 
Z 
\subset 
\Bbb 
G(M)$ 
is 
the 
zero-set 
of a 
$G^{\Bbb 
C}$-invariant\footnote{ Please note that $G^{\mathbb C}$ is unimodular.}
section $s \in \Gamma(M,L_M^{-m}\otimes \pi^*{\mathcal
O}(\ell))$.

\noindent With this set-up the proof of \cite{M2},Theorem 1, where an 
argument in
(1.2) using
duality between Hermitian symmetric manifolds of compact and noncompact
type was used, generalizes to complete the proof of Theorem 1.
The starting point is the Poincar\'e-Lelong equation on the Hermitian
symmetric manifold of compact type $M$.  Letting $g_c$ be a canonical 
K\"ahler-Einstein metric on $M$, $\widehat{g}_c$ the induced Hermitian
metric on the tautological 
line 
bundle, 
$h_c$ 
be a 
$G_c$-invariant 
Hermitian
metric on $\mathcal 
O(1)$, 
we 
have, by 
the 
Poincar\'e-Lelong 
equation

$$
\frac{\sqrt{-1}}{2\pi} \partial\overline{\partial}\log\|s\|_c^2
= mc_1(L,\widehat g_c) - 
\ell c_1(\pi^*\mathcal 
O(1), 
\pi^*h_c) 
+ 
[\mathcal 
Z_M].
$$

\noindent (There was a sign mistake on the curvature term in the 
Poincar\'e-Lelong
equation in [\cite{M2},(1.2) and the proof of Theorem 1] but they do not affect
the rest of the argument.)  Let now $\widehat N \subset \Bbb G(M)$ be 
the canonical lifting of 
the 
totally-geodesic 
submanifold 
$N 
\subset 
M$ to 
the 
Grassmann bundle $\Bbb G(M)$. Let $\omega_c$ be the K\"ahler form of 
$(M,g_c)$. 
Then, from $\widehat N \cap \mathcal 
Z_M = 
\emptyset$ 
and 
Stokes' 
Theorem
we have

$$
\int_{\widehat N} \Big(mc_1(L,\widehat g_c) - 
\ell c_1(\pi^*\mathcal 
O(1), 
\pi^*h_c)\Big) 
\wedge 
(\pi^*\omega_c)^{p-1} = 
0,
$$

\noindent Denote by $\sigma: N \to \Bbb G(M)$ the canonical
lifting 
map 
whose 
image is precisely $\widehat N$.  We note that $\sigma^*(L,\hat g_c)$
is nothing other than $(K_{N},det(g_c|_{N}))$.  In fact there is 
a correspondence between the vector bundle 
$\Lambda^pT_N$ over $N$ and the tautological 
line 
bundle 
$L$ 
over 
$\Bbb P(\Lambda^pT_N)$, which sends a vector $\eta \in \Lambda^pT_x(N)$
to the vector, also denoted by $\eta$, as an element over the 
point $[\eta]$.  In this tautological 
identification 
the 
length of
vectors is preserved, so that $\sigma^*(L,\hat g_c)$
is nothing other than $(K^{-1}_{S_c},det(g|_{S_c}))$, as asserted. 
$\omega_c$ is a positive multiple of $c_1(\mathcal
O(1),h_c)$ 
and 
we 
may 
take the two to be the same in what follows.  Then,

\begin{eqnarray*}
0&=&\int_{\widehat N} \Big(mc_1(L,\widehat g_c) -
\ell c_1(\mathcal O(1), 
\pi^*h_c)\Big) 
\wedge 
(\pi^*\omega_c)^{p-1} \\
& =& \int_{N} \Big(mc_1(K^{-1}_{N},det(g_c|_{N})) -
\ell c_1(\mathcal 
O(1),h_c)\Big) 
\wedge 
\omega_c^{p-1}\\
& = &\int_{N} mRic(g|_{N}) -
\ell c_1(\mathcal 
O(1),h_c)\Big) 
\wedge 
\omega_c^{p-1}\\
& = &\int_{N} \Big(\frac{\phantom{,}m\phantom{,}}{p}K(g_c|_{N}) -
\ell\Big) \omega_c^p,
\end{eqnarray*}

\noindent where $K(g_c|N)$ stands for the scalar
curvature, 
which is 
a 
constant, 
which
forces $K(g_c|_{N}) = \frac{p\ell}{m} $.
Denote by $g$ the 
K\"ahler-Einstein metric on $\Omega$ dual to $g_c$, etc., and $(E,h)$ the 
negative line bundle on $\Omega$ dual to $({\mathcal
O}(1),h_c)$.
By our choice of $h_c$ we have $c_1(E,h)=-\omega$.
For the same section $s \in \Gamma(\Bbb G(M),L \otimes \pi^*E)$, restricted to 
$\Omega$, denote by $\|s\|$ the norm measured in terms of $g$ and 
$h$. Then, we have the Poincar\'e Lelong equation on $\Omega$
$$
\frac{\sqrt{-1}}{2\pi} \partial\overline{\partial}\log\|s\|^2
= mc_1(L,\widehat g) - \ell c_1(\pi^* E, \pi^*h) + [{\mathcal
Z}_{\Omega}].
$$

\noindent Since    $s$ is $G^{\mathbb C}$ invariant, every
ingredient of this Poincar\'e Lelong equation descends to
 ${\mathbb G}(\Omega) / \Gamma$ and gives rise to a relation between 
cohomology classes
 on this manifold:

 $$ [{\mathcal
Z}_{\Omega /\Gamma}]=-mc_1(L)_{\Omega /\Gamma} + \ell c_1(\pi^* E_{\Omega 
/\Gamma}) .$$

\noindent Suppose now we have a compact totally-geodesic complex submanifold
$S_o \subset
X := \Omega/\Gamma$ uniformized by $D \subset \Omega$. Then, $\widehat S_o 
\cap {\mathcal Z}_X = \phi$
and we have

\begin{eqnarray*}
0 & = &\int_{\widehat S_o} \Big(mc_1(L,\widehat g) -
\ell c_1(E, \pi^*h)\Big) \wedge (\pi^*\omega)^{p-1} \\
& = &\int_{S_o} \Big(mc_1(K^{-1}_{S_o},det(g|_{S_o})) -
\ell c_1(E,h)\Big) \wedge \omega^{p-1}\\
& = &\int_{S_o} mRic(g|_{S_o}) - \ell c_1(E,h)\Big) \wedge \omega^{p-1}\\
& = &\int_{S_o} \Big(\frac{\phantom{,}m\phantom{,}}{p}K(g|_{S_o}) +
\ell\Big) \omega^p.
\end{eqnarray*}

\noindent Since $K(g|_{S_o}) = -K(g_c|_N)$ the integrand vanishes identically
on 
$S_o$, which has to be the case {\it a priori}.  But we are going 
to characterize $p$-dimensional compact complex submanifolds $S \subset X$
whose tangent spaces do not belong to $\mathcal 
Z_X$.  
For 
the 
proof we {\it do not} assume the existence of $S_o$.  Integrating the 
restriction of the Poincar\'e-Lelong equation to $S$ and assuming
that $\widehat{S}_o \cap {\mathcal Z}_X = \phi$ we conclude that
$$      
\int_{S} \Big(\frac{\phantom{,}m\phantom{,}}{p}K(g|_{S}) +
\ell\Big) \omega^p = 0
$$

\noindent Once we have the intepretation of the integral over $\widehat{S}$ as
an integral over $S$, Theorem 3 follows from \cite{Eys1}, Proposition 9.2.5.
 The rest of the argument is actually a standard generalization
of the proof of the Arakelov inequality, which
we
include for the sake of completeness.
Indeed,  the scalar curvature of $D$
is $K(g|_D)=-C \|\Sigma(T_o D) \|^2$ by standard curvature formulas
(see e.g. \cite{Eys1} p. 205-206), $C$ being a positive constant,
and the scalar curvature of $S$ at $x$
is $K(g|_S)_x= -C\|\Sigma (T_x S)\|^2 - \| \sigma _x\|^2$ where $\sigma$ is the
second fundamental form  of the embedding $S\subset \Omega/ \Gamma$.
By Lemma \ref{absmin} we have, at any point $x$,
 $\frac{\phantom{,}m\phantom{,}}{p}K(g|_{S})_x +
\ell  \le -\frac{\phantom{,}m\phantom{,}}{p} \| \sigma _x\|^2$.
Since the  integral of the l.h.s. is $0$, we deduce that $\sigma_x=0$,
i.e. $S$ is totally geodesic. This concludes the proof of Theorem 3. \quad
$\square$

\noindent
{\bf Remarks:}

\rm\noindent
(a) It was convenient to make use of the compact dual and duality to check
that for $s \in \Gamma(\Bbb P(T_{\Omega}), L^{-m} \otimes \pi^*(E))$,
the ratio of $m$ to $\ell$ is the right one.  The same thing can
be obtained by working on $\Omega$ alone, provided that we assume the 
 fact that $D$ 
admits 
a
torsion-free 
cocompact
lattice, a result
dating back to Borel (\cite{Bo}, Corollary to Theorem A).

\noindent (b) Modifying the above proof with a brute force curvature
computation of the appropriate line bundle would also allow us to 
prove Theorem 3 for reducible domains $\Omega$. We did not feel compelled
 to give the full argument for such a slight generalization 
 we nevertheless have to mention.

\subsection{Explicit examples}

\noindent In this section, we give some examples of pairs $(\Omega,D)$ of 
bounded symmetric
domains, $\Omega$ irreducible, $\Omega = G/K$, $\dim(D) = p$, such that 
$[T_o(D)]$ is a semistable point in the Grassmannian 
$\text{{\rm Gr}}(p,T_o(\Omega))$ of $p$-planes in $T_o(\Omega)$.  By Theorem 
3, gap
rigidity in the Zariski topology holds for such pairs $(\Omega,D)$. In fact,
we do more: we construct an explicit invariant hypersurface. 

\paragraph{(1) Variation of Hodge Structures and Gauss-Manin complexes}
Let $\Omega$ be a bounded symmetric domain and $D \subset \Omega$
be a 1-hyperrigid domain in the sense of Eyssidieux \cite{Eys2}
from variation of Hodge structures.  The methods and results there apply 
more generally to period domains.  \cite{Eys2}
contains tables of lists
of 1-hyperigid subdomains in the case of bounded symmetric domains, but not a 
full classification.
By \cite{Eys1},
Proposition 9.3.6, gap rigidity holds for $(\Omega,D)$ in the Zariski
topology, which was proven there by means of the Gauss-Manin complex.
There an excluded hypersurface ${\mathcal Z}_o\subset
\text{{\rm Gr}}(p,T_o(\Omega))$ 
such that $[T_o(D)] \notin {\mathcal
Z}_o$ 
can 
be 
identified 
in 
terms of
Lie 
algebras.
Alternatively, for the Gauss-Manin complex $(X,K_p^*)$ arising from the 
variation of Hodge structures, we can find an explicit invariant hypersurface
 ${\mathcal
Z}_o$ 
of 
$p$-dimensional 
vector
subspaces in $\text{{\rm Gr}}(p,T_o(\Omega)$: this is the locus over which the 
Gauss-Manin complex
fails to be an exact sequence at the point $x \in X$
be described as the zero set at $x$.
 In other words this the zero set of the determinant of the 
complex\footnote{This
  is a invariant holomorphic section on $
{\rm Gr}(p, T_o\omega)$ of the homogenous line bundle defined by the 
determinant 
of the
cohomology, see \cite{GKZ}
for a beautiful exposition of this classical construction of Cayley. }.
The condition on scalar
curvatures
is
verified
in
the 
calculation
there.

\paragraph{(2) Holomorphic curves on certain irreducible bounded symmetric
domains}
The set of all irreducible bounded
symmetric domains $\Omega$ for which there is a $K$-invariant hypersurface
is listed in \cite{M2}, Proposition 1.  For dim$(D) = 1$, gap rigidity
holds
for $(\Omega,D)$ in the Zariski topology if and only if $\Omega$ belongs
to that list and $D$ is a totally-geodesic disk of maximal type, i.e., 
equivalent under $\text{Aut}(\Omega)$ to the diagonal disk of a maximal 
polydisk
(which is of dimension equal to the rank of $\Omega$).  To put this in the 
framework of Theorem 3 we list here such domains $\Omega$ together with some
invariants.

\noindent They will be used again in (4). In the following we normalize the
K\"ahler-Einstein curvature so that the minimal disk is of constant 
Gaussian curvature $-2$.  With this normalization the K\"ahler-Einstein 
constant $\rho_{\Omega}$ agrees with $-c_1(M)$, where $c_1(M)$ is the 
first Chern class of the compact dual $M$, identified in the standard
way with a positive integer.

\begin{description}
\item{\rm (a)} 
$\Omega$ of Type $I_{m,n}$ with $m = n > 1;  \ r = n,
\ \text{dim}(\Omega) = n^2, \ \rho_{\Omega} = -2n, $    \\
\item{\rm (b)}
$\Omega$ of Type $II_n$ with $n$ even; $\ r = 
\frac{\phantom{,}n\phantom{,}}{2}, \ 
\text{dim}(\Omega) = \frac{n(n-1)}{2},
\ \rho_{\Omega} = -2(n-1)$;    \\
\item{\rm (c)}
$\Omega$ of type $III_n, n \geq 2$; $\ r = n, \ \text{dim}(\Omega) = 
\frac{n(n+1)}{2}, 
\ \rho_{\Omega} = -(n+1)$; \\
\item{\rm (d)}
$\Omega$ of Type $IV_n, n \geq 3$; $\ r = 2, \ \text{dim}(\Omega) = n, \ 
\rho_{\Omega} = -n;$
and\\
\item{\rm (e)}
$\Omega$ of Type $VI$ (the 27-dimensional exceptional domain pertaining to
$E_7$); $\ r = 3, \ \text{dim}(\Omega) = 27, \ \rho_{\Omega} = -18.$
\end{description}

\noindent Note that the domains $\Omega$ are precisely those of tube type.
In retrospect, this can be conceptually explained by the fact that a domain is 
of tube type iff
the diagonal disk embedding is $(H_2)$, \cite{Sat}, p. 150, Remark 1.

\paragraph{(3) Holomorphic quadric structures}

\noindent
Let $n \geq 3$ and $\Omega = D^{IV}_n$ denote the $n$-dimensional bounded 
symmetric domain of 
Type IV, which is dual to the $n$-dimensional hyperquadric.
In Mok \cite{M2} \S4, Theorem 4,  we proved that gap rigidity holds for
$(D^{IV}_n,D^{IV}_p,i)$
in the Zariski topology for $1 \leq p < n$ and for $i: D^{IV}_p \hookrightarrow
D^{IV}_n$ the standard embedding.  For $p = 1$ this was included in (2).  
For $p > 1$ we exploited the existence and uniqueness of K\"ahler-Einstein
metrics on projective manifolds with ample canonical 
line 
bundles 
and 
made 
use
of holomorphic quadric structures. In this case $[T_o(D)] \in 
\text{{\rm Gr}}(p,T_o(\Omega))$, 
$D = D^{IV}_p$ is a 
semistable point under $K$-action in view of the canonical 
holomorphic 
quadric 
structure $Q$ on $\Omega$.  Here at each $x \in \Omega$ we have 
$Q_x : S^2T_x \to E_x$ where $E$ is homogeneous holomorphic line bundle on 
$\Omega$.  The excluded subvariety ${\mathcal Z}_o\subset
\text{{\rm Gr}}(p,T_o(\Omega)$ 
consists precisely of those $p$-planes $V$ such that 
$Q_o|_V$ is nondegenerate.  ${\mathcal 
Z}_o$
is 
the 
zero 
locus 
of a 
holomorphic 
section of a positive power of the dual tautological 
line 
bundle 
on 
$\text{{\rm Gr}}(p,T_o(\Omega))$ corresponding to the discriminant of the 
(twisted) 
complex bilinear form $Q_x$, which is clearly a hypersurface.  

\paragraph{(4) Maximal polydisks on certain irreducible bounded symmetric 
domains}

\noindent
Let $\Omega$ be one of the irreducible bounded symmetric domains of 
characteristic codimension 1 as listed in \cite{M2}, Proposition 1.
Let  $r > 1$ be its rank.
Denote by $\Delta^r$ a maximal (totally-geodesic) polydisk in $\Omega$.
By (2), gap rigidity holds in the Zariski topology for  
$(\Omega,\delta(\Delta^r))$, where $\delta(\Delta^r)$ stands for the 
diagonal disk.  Actually, gap rigidity also holds for $(\Omega,\Delta^r)$,
giving new examples to which Theorem 3 applies.  To see this,
 as to semistability of $[T_o(\Delta^r)]$, by \cite{M2}, Proposition 3, there
exists a
$G$-invariant holomorphic section $s$ of $L^{-r} \otimes \pi^{*}E$ over
$\Bbb P(T_{\Omega})$, where $L$ stands for the tautological
line
bundle 
over
$\Bbb P(T_{\Omega})$, $\pi: \Bbb P(T_{\Omega})\to \Omega$ denotes the 
canonical projection, and $E$ is a 
homogeneous 
positive 
holomorphic 
line bundle over $\Omega$, such that the zero set of $s$ is precisely
the highest characteristic bundle $\pi: {\mathcal
S }
\to 
\Omega$.  
(Here 
there 
is 
a Zariski open $K^{\Bbb C}$ orbit ${\mathcal
O}_o$ 
in 
$\Bbb 
PT_o(\Omega)$ 
and 
${\mathcal S}_o \subset \Bbb P(T_o(\Omega))$ is its complement.)  
Equivalently 
$s$ corresponds to a $G$-invariant holomorphic section $\sigma$ of the 
homogeneous holomorphic vector bundle $S^r T_{\Omega}^* \otimes E$ on $\Omega$.
Writing $T$ for $T_{\Omega}$, $\sigma$ induces a $K$-invariant linear map 
$\theta: T \to S^{r-1}T^* \otimes E$, and hence a $K$-invariant map
$\wedge^r\theta: \Lambda^rT \to \Lambda^r(S^{r-1}T^*)\otimes E^r$. For a 
complex vector space $V$ and positive intergers $n$ and $m$ let 
$\mu: \Lambda^n(\otimes^mV) \to \otimes^m(\Lambda^nV)$ be the linear 
map defined by

\begin{eqnarray*}
& \qquad \mu\big((v_{11}\otimes\cdots\otimes v_{1m}) \wedge \cdots \wedge
(v_{n1}\otimes\cdots\otimes v_{nm})\big)\\
& = (v_{11}\wedge\cdots\wedge v_{n1})\otimes\cdots\otimes 
(v_{1m}\wedge\cdots\wedge v_{nm}).
\end{eqnarray*}

\noindent If $v_{ij}$ is independent of $j$ then the image under $\mu$ lies
in $S^m(\Lambda^nV)$.  By polarization we conclude that by
restriction $\mu: \Lambda^n(S^mV) \to S^m(\Lambda^n V)$.  
Denoting also by $\mu$ the corresponding bundle homomorphism 
applied to $T$ we obtain a homomorphism 
$\nu = (\mu\otimes id_{E^r})\circ \wedge^r\theta$,
$\nu: \Lambda^rT \to S^{r-1}(\Lambda^rT^*)\otimes E^r$, which 
therefore defines a $K$-invariant element of 
$\Lambda^rT^*\otimes S^{r-1}(\Lambda^rT^*)\otimes E^r$, and 
hence by canonical 
projection 
a 
$K$-invariant 
element 
$\tau \in S^r(\Lambda^rT^*)\otimes E^r$.  Denoting also by 
$\pi: \Bbb G(\Omega) \to \Omega$ the Grassmann bundle whose
fiber over $x \in \Omega$ is the Grassmannian of $r$-planes
in $T_x(\Omega)$ we have obtained a $G$-invariant holomorphic
section $t \in \Gamma(\Bbb G,L^{-r}\otimes \pi^*E^r)$,
where $L$ denotes here the tautological 
line 
bundle 
over
$\Bbb G$.

\noindent It remains to show that $t([T_o(\Delta^r)]) \notin 0$.
Given this, the zero set of $\tau$, which is necessarily
non-empty and of codimension 1 in $\Bbb G_o = \text{{\rm Gr}}(r,T_o(\Omega))$,
defines the excluded hypersurface ${\mathcal
Z}_o$, 
showing 
that 
gap rigidity holds in the Zariski topology for $(\Omega,\Delta^r)$.  
To this end 
for notational convenience consider the case where
$\Omega$ is a Type I domain $D^I_{r,r}$ with $r > 1$, so 
that the tangent space can be identified as the space of 
$r$-by-$r$ matrices. We note that the same argument works in 
general for any $\Omega$ of characteristic codimension 1.
Denote by $e_{ij}$ the $r$-by-$r$ matrix whose entries are 0
except for the $(i,j)$-th entry, which is equal to 1.  
Thus $(e_{ij})$ constitutes a basis for $T_o(\Omega)$, whose 
dual basis will be denoted by $(e^*_{ij})$.  For 
the maximal polydisk $\Delta^r$ we may take $T_o(\Delta^r)$
to consist precisely of the diagonal matrices.  Identifying 
$E_o$ with $\Bbb C$ the $K$-invariant section $\sigma$ corresponds
to the determinant function. Thus $\theta: T \to S^{r-1}T^* \otimes E$
satisfies $\theta(e_{ii}) = e_{11}^*\circ\cdots \circ\widehat{e^*_{ii}}
\circ\cdots\circ e^*_{rr} +$ terms vanishing on $T_o(\Delta^r)$,
 where $\circ$ denotes the symmetric product
and\quad $\widehat\cdot$\quad denotes exclusion. Hence, 

\begin{eqnarray*}
& \wedge^r\theta(e_{11}\wedge\cdots\wedge e_{rr}) = (e^*_{22}\otimes  e^*_{33}
\otimes \cdots \otimes e^*_{rr})\wedge (e^*_{11}\otimes  
e^*_{33}\otimes \cdots \otimes e^*_{rr})\\
& \qquad\qquad \wedge \cdots \wedge 
(e^*_{11}\otimes e^*_{22}\otimes \cdots \otimes e^*_{r-1,r-1})+
\ {\text{terms vanishing on}}\  T_o(\Delta^r)
\end{eqnarray*}

\noindent It follows by a straightforward calculation
that
$\nu(e_{11}\wedge\cdots\wedge e_{rr})$ is a positive multiple
of $(e^*_{11}\wedge \cdots \wedge e^*_{rr})^{r-1}$ modulo
terms vanishing on $T_o(\Delta^r)$, implying that 
$t(e_{11}\wedge\cdots\wedge e_{rr})^r) \neq 0$ at $[T_o(\Delta^r)])$, as 
desired.

\noindent If $r = k\ell$ for positive integers $k, \ell > 1$, we can write
the maximal polydisk as $\Delta^r = (\Delta^k)^{\ell}$. Then, the same
argument as in the above shows that gap rigidity holds in the Zariski topology 
for the pair $\big(\Omega, \big(\delta(\Delta^k\big)^{\ell}\big)$ where
$\delta(\Delta^k) \cong \Delta$ is the diagonal.

\paragraph{(5) Examples where the subdomains are higher-dimensional and
irreducible}

\noindent
The arguments of (4) give rise to examples of pairs $(\Omega,D)$ with 
$\Omega = D^{I}_{r,r}$, where $D$ is higher-dimensional and irreducible.  
These are the subdomains $D = D^{II}_r$, consisting of skew-symmetric 
matrices, and the subdomains $D = D^{III}_r$, consisting of 
symmetric matrices\footnote{In some cases, $D \subset \Omega$ occurred
as a 1-hyperrigid
subdomain, namely $D^{II}_r$ and $D^{III}_r, \ r\cong 0,1 [4]$, see 
\cite{Eys2}.}.
 We have $\dim(D^{II}_r) = \frac{r(r-1)}{2}$ and
$\dim(D^{III}_r) = \frac{r(r+1)}{2}$.   As in (4) in the process the excluded
hypersurface
${\mathcal Z}_o \subset \text{{\rm Gr}}(p,T_o(\Omega))$ can in principle 
be 
explicitly 
determined.  
For the purpose of illustration we will 
establish the semistability of $[T_o(D)]$ in $\text{{\rm Gr}}(p,T_o(\Omega)$ 
in 
the case of $r = 3$.

\noindent We start with $D = D^{III}_3$, which is a 6-dimensional bounded 
symmetric
domain of rank 2,  Let $(e_{ij})$ be a basis of $T_o(\Omega)$ and $(e^*_{ij})$ 
be a dual 
basis, as in (4).  Then, $T_o(D)$ is spanned by the basis 
$\big\{e_{12}+e_{21}$, $e_{13}+e_{31}$, $e_{23}+e_{32}$, 
$e_{11},e_{22},e_{33}\big\}$.
We write $x_{ij}$ for $e_{ij}$ mod $T_o(D)$.  For the determinant $det$ on 
3-by-3 matrices, we have
$$
det = e_{11}^*e_{22}^*e_{33}^* + e_{12}^*e_{23}^*e_{31}^*
+ e_{13}^*e_{21}^*e_{32}^* - e_{13}^*e_{22}^*e_{31}^*
- e_{11}^*e_{23}^*e_{32}^* - e_{12}^*e_{21}^*e_{33}^*, 
$$
so that
$$
\theta(e_{12}) = e_{23}^*e_{31}^* - e_{21}^*e_{33}^*, \quad
\theta(e_{21}) = e_{13}^*e_{32}^* - e_{12}^*e_{33}^*.
$$

\noindent Denote by $\overline{\theta}$ the composite $\kappa \circ \theta$,
where $\kappa$ is the projection map induced by the 
quotient homomorphism $T^*_o(\Omega) \to T^*_o(D)$.  
We have

$$
\overline{\theta}(e_{12}) = x_{23}x_{31} - x_{21}x_{33},  \quad 
\overline{\theta}(e_{21}) = x_{13}x_{32} - x_{12}x_{33}.
$$
\noindent Noting that $x_{ij} = x_{ji}$ we have
$$
\overline{\theta}(e_{12}) = x_{23}x_{13} - x_{12}x_{33} =
\overline{\theta}(e_{21})\ ,
$$
\noindent implying by analogous calculations

\begin{eqnarray*}
\overline{\theta}(e_{12}+e_{21}) &=& 2(x_{23}x_{13} - x_{12}x_{33}); \\
\overline{\theta}(e_{13}+e_{31}) &=& 2(x_{12}x_{23} - x_{13}x_{22}); \\
\overline{\theta}(e_{23}+e_{32}) &=& 2(x_{12}x_{13} - x_{23}x_{11})
\end{eqnarray*}

\noindent On the other hand
$$
\theta(e_{11}) = e_{22}^*e_{33}^* - e_{23}^*e_{32}^*,
$$
\noindent implying by analogous calculations

$$
\overline{\theta}(e_{11}) = x_{22}x_{33} - x_{23}^2, \quad \quad
\overline{\theta}(e_{22}) = x_{11}x_{33} - x_{13}^2, \quad \quad
\overline{\theta}(e_{33}) = x_{11}x_{22} - x_{12}^2. 
$$

\noindent Putting together the formulae for $\overline{\theta}$ for the basis 
vectors
$\big\{e_{12}+e_{21}$, $e_{13}+e_{31}$, $e_{23}+e_{32}$, $e_{11}$, $e_{22}$,
$e_{33}\big\}$, 
we obtain a formula for the section 
$t \in \Gamma(\Bbb G_o,L^{-3})$. Denoting by $\overline{\mu}$
the homomorphim induced by $\mu: \Lambda^6(S^2T_o) \to S^2(\Lambda^6T_o)$
as in (4), and by the quotient homomorphism $T_o^* \to T_o^*(D)$, we have

\begin{eqnarray*}
& \qquad t\big((e_{12}+e_{21})\wedge (e_{13}+e_{31})\wedge (e_{23}+e_{32}) 
\wedge e_{11}
\wedge e_{22} \wedge e_{33}\big) \\
& = 8\overline{\mu}
\big((x_{23}x_{13}-x_{12}x_{33}) \wedge (x_{12}x_{23}-x_{13}x_{22})\wedge 
(x_{12}x_{13}-x_{23}x_{11}) \\
& \qquad \wedge (x_{22}x_{33}-x_{23}^2)
\wedge (x_{11}x_{33}-x_{13}^2) \wedge (x_{11}x_{22}-x_{12}^2)\big).
\end{eqnarray*}

\noindent The only terms that count come from the combination of the first
halves of each of the 6 quadratic expressions inside brackets, since all other
exterior products are zero, as the variables $(x_{23}, x_{12}, x_{13},
x_{11}, x_{22}, x_{33})$ do not occur the right number of times.  We get
finally

\begin{eqnarray*}
& \qquad
t\big((e_{12}+e_{21})\wedge (e_{13}+e_{31})\wedge (e_{23}+e_{32}) \wedge e_{11}
\wedge e_{22} \wedge e_{33}\big)\\
& = 32(x_{23}\wedge x_{12} \wedge x_{13}
\wedge x_{22} \wedge x_{33} \wedge x_{11})^2 \neq 0,
\end{eqnarray*}

\noindent as desired.  The zero set of $t$ defines the cubic hypersurface
${\mathcal Z}_o \subset \Bbb G_o$ on the Grassmannian $\Bbb G_o$ 
which 
serves as 
the 
excluded subvariety in Theorem 3.

\noindent Next we consider $D = D^{II}_3$, which is a 3-dimensional bounded 
symmetric
domain of rank 1,  i.e., biholomorphic to the 3-dimensional unit ball $B^3$.
A straightforward adaptation of the preceding argument does not work. 
Use the same notations as in the above, except that $D$ stands for 
$D^{II}_3$ and that $x_{ij}$ is $e_{ij}^*$ modulo $T_o(D)$, i.e., 
the vector subspace of skew-symmetric matrices in the vector space of 3-by-3 
matrices, 
so that $x_{ij} = -x_{ji}$, in particular $x_{33} = 0$, and we have

\begin{eqnarray*}
\overline{\theta}(e_{12}) = x_{23}x_{31}, \
\overline{\theta}(e_{21}) = x_{13}x_{32} \\
\overline{\theta}(e_{21}) = (-x_{31})(-x_{23}) = x_{23}x_{31} 
= \overline{\theta}(e_{12}), \quad {\text{hence}}\\
\overline{\theta}(e_{12} - e_{21}) = 0.
\end{eqnarray*}

\noindent It follows that we get the section $t \equiv 0$ in this case.
To overcome the difficulty, we work with symmetric matrices
in place of skew-symmetric matrices by going to the cotangent
bundle, as follows.  The starting point of the construction of 
$\cal Z$ is the existence of the determinant.  On the compact
dual $M$, which is here the Grassmannian $\text{{\rm Gr}}(r,V)$ of 
3-planes on $V \cong \Bbb C^6$, and at a point
$x \in M$, identifying $T_x(M)$ with the set of 
3-by-3 matrices, we can define tentatively the 
`determinant' of the tangent vector. However, the matrix reprentation 
is unique only up to the complexification $K^{\Bbb C}$
of the isotropy subgroup
$K$ at $o$.  The action of the centre of $K_x$ on the `determinant'
shows that there is a well-defined determinant on the projectivized tangent 
bundle of $M$, as an 
$\text{Aut}(M)$-invariant section of a homogenenous holomorphic 
line bundle.  Likewise we can define a determinant for the 
cotangent bundle as a section on $\Bbb P(T^*(M))$, with 
tautological line bundle 
$\Lambda$, of 
a 
homogeneous
holomorphic line bundle which restricts to $\Lambda^{-r}$ 
over $\Bbb P(T^*_o(M))$, corresponding to an element 
of $S^rT_o(M)$.  Now in general there is a canonical 
correspondence
between the Grassmann bundle $\Bbb G(M)$ over $M$ of 
$p$-dimensional tangent planes over $M$ with the Grassmann bundle
$\Bbb G'$ of $(m-p)$-dimensional vector subspaces of cotangent spaces over 
$M$, obtained by sending each $p$-plane $A \subset
T_x(M)$ to its annihilator $A^{\perp} \subset T^*_x$.  
Consider a reference point $o \in \Omega$ and represent
$T_o(\Omega)$ as the space of 3-by-3 matrices such that 
$T_o(D)$ is identified with the 3-dimensional vector subspace
of skew-symmetric matrices.  Identifying also $T^*_o(\Omega)$,
via the natural complex bilinear pairing between 
$T_o(\Omega)$ and $T_o^*(\Omega)$, as the vector space of 
3-by-3 matrices.  Then $(T_o(D))^{\perp}$ is nothing 
other than the 6-dimensional vector space of symmetric 
matrices.  Working with the determinant on the cotangent 
bundle in place of the tangent bundle we can perform the 
same proof as in the above to find a cubic hypersurface in 
$\Bbb G'_o$ which avoids the point $[T_o^{\perp}(D)]$.  The 
canonical isomorphism $\Bbb G_o(D) \cong \Bbb G'_o(D)$
then 
gives a cubic hypersurface $\mathcal Z_o \subset \Bbb G'_o$ which avoids 
$[T_o(D)]$.

\noindent
{\bf Remarks}

\rm
\noindent
(a) Examples discuss in $(2) - (5)$ are all obtained on irreducible bounded 
symmetric
domains $\Omega$ of characteristic codimension 1, i.e., those admitting a 
$K$-invariant 
hypersurface.  This does not need to be the case for $(\Omega,D)$ to exhibit 
gap rigidity for some choice of $D$.   Such examples occur in the context of 
Gauss-Manin complexes.  For certain values of $p, n > 1$, and for 
$\Omega = D^{I}_{n,pn}$ it was shown in Eyssidieux \cite{Eys2} that
$(D^{I}_{n,pn},B^p;i)$ is 1-hyperrigid, where $i: B^p \to D^I_{n,pn}$ is 
obtained
by embedding $B^p$ as the diagonal of a product $P = B^p \times \cdots \times 
B^p$ ($n$ 
factors) of complex unit $p$-balls, and $P$ is realized as a totally geodesic  
     submanifold of $\Omega$ in a standard way.    
 This  was precisely the hint towards the general theory.

(b) The method used in (4) and (5) can be applied to yield various examples 
where $(\Omega,D)$ exhibits gap rigidity in the Zariski topology.  They are
typically of `diagonal' type, for instance, if $n 
= 
p\ell$, 
then 
$\big(D^{III}_n,\delta\big((D^{III}_p)^{\ell}\big)\big)$ can be shown to 
exhibit
gap  rigidity in the Zariski topology.  The methods of (4) and (5) produce
in each case a locally 
homogenenous 
holomorphic 
section 
$t$ in 
some 
$\Gamma(\Bbb P(T_X),L^{-m}\otimes \pi^*E)$, and the question is to know whether
$t([T_o(D)]) \neq 0$.  As is shown in (5), it can happen that $t \equiv 0$,
although in that particular example there is a way to circumvent the 
difficulty. 
When however $t \not\equiv 0$, then there is the advantage that $t$ can in 
principle
be explicitly determined.  Moreover, by \cite{M2}, Proposition 3, whenever the
higher
characteristic subvariety $\mathcal 
S_o 
\subset 
\Bbb 
PT_o(\Omega)$ 
 is 
a 
hypersurface,
it is always of degree $r = rank(\Omega)$.  When $t\ne 0$, the construction in 
(4) and (5) 
yields a symmetric polynomial of degree $r$. Thus, in the case
of (3) on holomorphic quadric structures, the excluded hypersurface
${\mathcal Z}_o \subset {\mathbb G}_o$ is always the zero set of a
nontrivial 
quadratic
polynomial.

\subsection{Classification of $(H_3)$-embeddings into an irreducible domain.}

\noindent
An embedding into a
reducible domain is $(H_2)$ if and only if its factors are also 
$(H_2)$-embeddings. This reduces the classification of $(H_2)$-embeddings
to the classification of 
maximal  $(H_2)$-embeddings into an irreducible domain.
We  give two tables, extracted from \cite{Iha} and \cite{Sat}. 
The first one gives  all
maximal  $(H_2)$-embeddings into 
classical domains up to equivalence.
The second table gives
all maximal  $(H_2)$-embeddings into an exceptional
domains
 and, for every irreducible $(H_2)$-subdomain, a chain
relating it to a maximal one.\\

\noindent
$(H_3)$-embeddings $(\Omega, D)$ into an irreducible domain $\Omega$
are $(H_2)$ embeddings with an additionnal requirement on the 
Einstein constants of the  factors  of $D$ if $D$ is not irreducible. 
It is  straightforward to use the tables to give a complete list
of all $(H_2)$ embeddings into a given irreducible domain, it
 is easy in each case to
 compute the Einstein constant
of the induced metric on every irreducible component
and decide whether they agree, giving a complete list 
of $(H_3)$-embeddings into $\Omega$. Unfortunately, we had not enough time
to perform the required curvature computations. 
We will be content with the remark that
a way to check without calculation that a $(H_2)$-embedding of a reducible 
domain of the form $D=\Omega_1^{n}$, where $\Omega_1$ is irreducible, is $(H_3)$ when
the automorphism group of $\Omega$ permutes the various subdomains of the
form $o\times\ldots\times o \times \Omega_1 \times o \times \ldots$.
This gives $(H_3)$-embeddings of the form $(D^I_{kp,kq}, (D^I_{p,q})^k)$, $(D^A_{nk}, (D^A_n)^k)$
$A= II, III$ , $(D^{III}_{20p}, (B^5)^p)$, etc. We do not know any example
which does not fit this pattern.
$(H_3)$-embeddings into a reducible domain are also in principle straightforward
to classify.

\begin{tabular}{|c|c|c|c|}
\hline
\multicolumn{4}{|c|}{Maximal $(H_2)$-subdomains}\\
\multicolumn{4}{|c|}{of a classical domain}\\
\hline
$\Omega $ &  $D $&  maximal &Additional conditions\\
\hline
$D^I_{p,q}$ & $D^I_{r,s}\times D^I_{p-r,q-s}$  &* &
$\frac{r}{s} = \frac{p}{q}$\\
& & & $(H_3)$ iff $p=r$ \\
\hline  & $D ^{II}_n$  & * & $p=q=n$ \\
\hline & $D ^{III}_n$  & * & $p=q=n$ \\
\hline
& $B^m$  &$m\not = 2r+1$& $p={m\choose {r-1}} ,  \
q=  {m\choose r}, r\in {\mathbb N}$\\
\hline
& $D^{IV}_{2l}$ & $l \equiv 0 [2]$& $p=q= 2^{l}, l\ge 3$ \\
\hline
& $D^{IV}_{2l-1}$  & & $p=q= 2^{l-1},  l\ge 3$ \\
\hline
$D^{II}_n$ & $D^I_{r,r}$&*& $n=2r$ \\
\hline
           & $D^{II}_r \times D^{II}_{n-r} $ &*& $n>r$\\
	   &&& $(H_3)$ iff $n=2r$\\
           \hline
& $B^m$ &*& $n={m+1\choose\frac{m+1}2 }, m \equiv 3[4]$\\
\hline
& $D^{IV}_{2l}$ &*& $n= 2^{l}, l\ge 3, l\equiv 3 [4]$ \\
\hline
& $D^{IV}_{2l-1}$  &*& $n= 2^{l-1},  l\ge 3, l\equiv 0,3 [4]$ \\
\hline
$D^{III}_n $& $ D^I_{r,r}$ &*& $n=2r$ \\
\hline  
   & $D^{III}_r\times D^{III}_{n-r} $ &*& $n>r$\\
   &&& $(H_3)$ iff  $n=2r$\\
   \hline
& $B^m$&*& $n={m+1\choose\frac{m+1}2 }, m \equiv 1[4]$\\
\hline
& $D^{IV}_{2l}$ &*& $p=q= 2^{l}, l\ge 3, l\equiv 1 [4]$ \\
\hline
& $D^{IV}_{2l-1}$  &*& $p=q= 2^{l-1},  l\ge 3, l\equiv 1,2[4]$ \\
\hline
$D^{IV}_{2l}$ & $D^I_{2,2}$& & $l\ge 3$\\
\hline 
              & $D^{IV}_{2l-1} $&*& $l\ge 3$ \\
\hline
$D^{IV}_{2l-1}$ & $D^{IV}_{2l-2}$ &*& $l\ge 3$\\
\hline
\end{tabular}\\

\pagebreak

\begin{tabular}{|c|c|c|c|}
\hline
\multicolumn{4}{|c|}{Maximal and irreducible
$(H_2)$-subdomains }\\
\multicolumn{4}{|c|}{of exceptionnal domains  }\\
\hline
$\Omega $&$ D $& $(H_3)$&  Chains of $(H_2)$-subdomains\\
\hline
 $D^V $& $D^I_{2,4}$ & *& $B^2 \subset B^2\times B^2  \subset D^I_{2,4}$\\
\hline
 &$ B^5 \times \Delta  $& &   \\
\hline
 $D^{VI}$ & $B^5 \times B^2$ & & \\
\hline
& $D^I_{2,6}$&  *& 
 $B^3 \subset B^3\times B^3 \subset  D^I_{2,6}$\\
\hline
 & $D^I_{3,3} $ &*&
 $\Delta \subset \Delta^3 \subset D^{III}_3 \subset D^I_{3,3}$  \\
\hline
 & $D^{II}_6 $ & *& $\Delta \subset \Delta^3 \subset D^{II}_6 $\\
\hline
 & $D^{IV}_{10}  \times \Delta $ &  &
$\Delta \subset \Delta^3 \subset D^{IV}_{10}  \times \Delta $\\
\hline
\end{tabular}

\section{ Overview on gap rigidity}

\subsection{Gap rigidity for $D^{III}_2$}
\noindent
 From  \cite{M2},Theorem 4,  (or Theorem 3 here) and Theorem 1,
 we have a complete understanding of gap rigidity for the unique 
3-dimensional bounded symmetric domain $\Omega$ of rank $> 1$. $\Omega$ can be 
described as a Type III domain of rank 2, equivalently the 3-dimensional 
Siegel 
upper half-plane, or as the 3-dimensional bounded symmetric domain of Type IV,
i.e., the noncompact dual of the 3-dimensional hyperquadric.

\begin{prop}\label{p1}
The question on the validity of gap rigidity on the 3-dimensional 
irreducible bounded symmetric domain $\Omega$ of rank $> 1$ is completely 
settled, 
as follows.  Denote by $\Delta^2 \subset \Omega$ a maximal totally-geodesic 
bidisk as given by
the Polydisk Theorem. There are, up to isometry, 
precisely 3 different types of positive-dimensional
totally-geodesic proper complex submanifolds $D$ of
$\Omega$, namely: 

\item{\rm (1)} $D = \Delta^2 \subset \Omega ;$

\item{\rm (2)} $D = \delta(\Delta) \subset \Delta^2 \subset \Omega ;$

\item{\rm (3)} $D = \Delta \times \{0\} \subset \Delta^2 \subset \Omega$.

\noindent
Gap rigidity holds in the Zariski sense  
for $(\Omega,D)$ for $D$ in {\rm (1)} or {\rm (2)}; 
but fails $($in the complex topology$)$ for $D$ in {\rm (3)}.
\end{prop}

\subsection{Rank one domains}

\noindent Theorem 3 and the examples in (2.5) show that gap rigidity in the
Zariski topology for a pair $(\Omega,D)$ can hold due to algebraic 
conditions satisfied by tangent planes to bounded symmetric domains.
Theorem 1 shows that gap rigidity can fail in the complex topology
due to product structures, although the construction of counterexamples
arising from holomorphic maps between compact Riemann surfaces of higher
genus does not generalize easily (cf. (3.3)).  There is a situation which
belongs
to neither of these situations and for which gap rigidity in the 
Zariski topology does not make sense.  This is especially the case
for bounded symmetric domains of rank 1.

\begin{ques}
Let $k < n$ be positive integers and embed the complex unit $k$-ball
$B^k$ into the complex unit $n$-ball $B^n$ in the standard way as a 
totally geodesic complex submanifold. Does gap rigidity 
hold for $(B^n,B^k)$ in the complex topology?
\end{ques}

\noindent There is up to this point no evidence as to whether one should
expect a positive or negative answer to Question 1.
Since the information of first order is trivial in this problem,  perhaps
one should try and
find a way to deal with higher order information,
which we can't do for the moment being.
The case
of $k = 1$ is perhaps the most difficult.  A negative
answer for $k = 1$, $n = 2$ constructed on quotients of the two-ball 
by torsion-free lattices would give first examples of exceptional divisors 
other than totally-geodesic cycles on some projective manifolds 
uniformized by the 2-ball.
 The case of $k > 1$ can be formulated
perhaps as a problem in two steps.  The first step is to understand 
whether a sufficiently pinched $k$-dimensional submanifold $S$ with 
$k > 1$ is necessarily uniformized by $B^k$.  This is equivalently
the question of asking whether $S$ admits a holomorphic projective 
structure, in view of \cite{KO}.  In the case of $k = 2$ it was
conjectured in Siu-Yang \cite{SY} that a compact K\"ahler-Einstein
surface with strictly negative sectional curvature is uniformized
by the complex unit 2-ball $B^2$.  In the same article, they
proved a type of pinching theorem which says that 
a compact K\"ahler-Einstein surface of nonpositive holomorphic bisectional 
curvature
must be biholomorphic to the 2-ball, provided that the curvature
at every point satisfies some very specific pinching condition.
It is possible to deform an $\epsilon$-pinched compact complex
submanifold $S \subset X := B^n/\Gamma$ to get a K\"ahler-Einstein 
surface with estimates on the deviation from constant holomorphic 
sectional curvature.  Unfortunately the pinching condition 
in \cite{SY} is not implied by such estimates.  The second step of the
problem is to prove that a holomorphic immersion of $S = B^k/\Gamma_o$
into $X = B^n/\Gamma$ is necessarily a totally geodesic immersion,
at least when the image is sufficiently pinched.  In the case where 
$2k > n$  this was settled in the positive in Cao-Mok \cite{CM} without
assuming any pinching condition.  In particular, in the case where 
$k = 2, n =3$ Question 1 is reduced to the first step.

\subsection{Embeddings into products}
\noindent
In relation to the counterexamples  showing that
gap rigidity does not hold for $(\Delta^2,\Delta \times o)$, one
can ask the question as to whether this type of construction can 
generalize.  Obviously one can get holomorphic mappings between 
products of compact Riemann surfaces of genus $> 2$ to show that 
gap rigidity fails for $(\Omega,\Delta^k;i)$ with rank$(\Omega)
= r > 1$, $\Omega$ not necessarily irreducible, $k < r$ and 
$i: \Delta \to \Omega$ arising from 
the inclusion $\Delta^k = \Delta^k\times \{0\} \subset \Delta^r$,
and an embedding of $\Delta^r$ as a maximal polydisk in $\Omega$.
On the other hand, for $\Omega$ an irreducible bounded symmetric
domain of rank $ > 1$, we have

\begin{prop}    \label{rignondiag}
Let $\Omega$ be an irreducible bounded symmetric domain of 
rank $> 1$.  Let $k < n$ be positive integers.  Then, gap 
rigidity holds for the pair $(\Omega^n,\delta(\Omega^k) \times \{0\})$  
in the Zariski topology when we consider only ambient complex 
manifolds $X$ of the form $\Omega^n/\Gamma$ with  
$\Gamma \subset Aut_o(\Omega^n)$. 
\end{prop}

\begin{prv}
Write $\Omega = G/K$ in the usual notations.
Here we can define, in the notations analogous to the statement 
of Theorem 3, ${\mathcal
Z}_o 
\subset 
\Bbb 
G_o$ 
to 
consist 
of 
tangent $p$-planes $V,p =\dim(\Omega)$, such that the 
canonical projection of $V$ onto, say 
the 
first 
factor 
$\Omega$
is an isomorphism.  With this convention, given any $S \subset 
X/\Gamma$, where $\Gamma \subset Aut_o(\Omega^n)$, the 
canonical projection $\rho: \Omega^n 
\to 
\Omega$ 
onto 
the 
first
factor induces on $S$ an integrable $K^{\Bbb C}$-structure.  It follows 
that the universal covering space of $S$ is biholomorphic to 
$\Omega$ (cf.  \cite{MY}). Thus, we have $S = \Omega/\Gamma_o$
as an abstract complex manifold, together with a holomorphic 
embedding $f: S \to X$. By Hermitian rigidity of Mok \cite{M1} it follows that 
$f:
S \to X$ 
is a totally geodesic isometric embedding up to a nonzero
normalizing constant.  Identifying $S$ with its image under
$f$ we conclude that $S \subset X$ is a holomorphic geodesic cycle.
\end{prv}

\noindent Proposition \ref{rignondiag} gives an example where gap rigidity can
hold in the
Zariski topology, with a reducible ambient domain $\Omega^n$,
even though the analogous condition on scalar 
curvatures 
as in 
Theorem 3 is {\it not\/} satisfied.  Here of course we are dealing with
reducible ambient domains, otherwise even the notion of gap 
rigidity in the Zariski topology as stated cannot be formulated,
but one may still raise the question, in the case of a pair 
$(\Omega,D)$ with $\Omega$  an irreducible bounded symmetric
domain, whether gap rigidity in the complex topology can 
hold for certain bounded symmetric domains $D$ which are say
irreducible and of rank $\ge 2$.  The problem can be reduced
to a question of holomorphic G-structures and is not within 
the scope of methods in the current article. 

\noindent Finally, it is intriguing even in the reducible case, whether
one can construct counterexamples to gap rigidity for pairs
$\big(B^n\times B^n,\delta(B^n\times B^n)\big)$ with $n > 1$.  We may 
formulate the problem as follows

\begin{ques}
Let $n > 1$.  Consider the set ${\mathcal
X}_n$ 
of 
all 
compact 
complex 
manifolds uniformized by the complex unit ball $B^n$.  Let
$\text{\rm Map}({\mathcal
X}_n)$ 
denote 
the
set 
of 
all 
nonconstant 
holomorphic 
mappings $f: X \to X'$ with $X, X' \in {\mathcal
X}_n$,
and 
$\text{\rm Map}_{\text{fin}}({\mathcal
X}_n)
\subset 
\text{\rm 
Map}({\mathcal
X}_n)$ 
the 
subset of 
all generically finite 
holomorphic 
maps. 
For each $f \in \text{\rm Map}({\mathcal
X}_n)$,
$f:X\to 
X'$, 
denote 
by 
$\mu(f) 
\in
(0,1]$
the real number defined by 
$\mu(f) = \text{\rm sup}\big\{\|df(x)\|: x \in X\big\}$.  
Does there exists a universal constant $c_n > 0$ depending only on $n$
such that $\mu(f) > c_n$ for any $f \in \text{\rm Map}_{\text{fin}}({\mathcal
X}_n)$
or more generally for $f \in \text{\rm Map}({{\mathcal
X}_n})$?
\end{ques}

\noindent
We note that $\mu(f) \leq 1$ because of the Ahlfors-Schwarz 
Lemma.  A negative answer to Question 2
for a given $n > 1$ will 
imply the failure of gap rigidity in the complex topology
for $(B^n\times B^n,B^n\times \{0\})$.  On the other hand,
a positive answer to an extended form of the question, 
where in place of $f \in \text{Map}({\mathcal
X}_n)$ 
we 
consider 
compact complex manifolds $X \in {\mathcal
X}_n$,
any 
representation
$\Phi: \pi_1(X) \to \text{Aut}(B^n)$ and any $\Phi$-equivariant
holomorphic map $F: B^n \to B^n$, will lead to the confirmation
of the gap rigidity in the complex topology for the 
pair $(B^n\times B^n,B^n\times \{0\})$.

\end{document}